\documentclass[a4paper,11pt]{article}       
\usepackage{graphics}
\usepackage{epsf}

\usepackage{amsfonts}
\usepackage{amssymb}
\usepackage{epsfig}
\usepackage[latin1]{inputenc}
\usepackage{color}
\usepackage{multirow}
\usepackage{array}
\usepackage[pagewise]{lineno}\linenumbers
\nolinenumbers
\usepackage{amsthm,amsmath,amssymb}

\newtheorem{theorem}{Theorem}

\providecommand{\keywords}[1]{\textbf{\textit{Keywords--}} #1}

\providecommand{\msc}[1]{\textbf{\textit{Mathematics Subject Classification--}} #1}

\date{}
\begin{document}

\title{Avoiding order reduction with explicit Runge-Kutta exponential methods in nonlinear initial boundary value problems}


\author{{\sc B. Cano} \thanks{Corresponding author. Email: bcano@uva.es}  \\ \small
IMUVA, Departamento de Matem\'atica Aplicada,\\ \small Facultad de
Ciencias, Universidad de
Valladolid,\\ \small Paseo de Bel\'en 7, 47011 Valladolid,\\ \small Spain \\
{\sc and}\\
{\sc M. J. Moreta}\thanks{Email:mjesusmoreta@ccee.ucm.es} \\
\small IMUVA, Departamento de An\'alisis Econ\'omico, \\
\small Facultad de Ciencias Econ\'omicas y Empresariales, Universidad
Complutense de Madrid, \\
\small Campus de Somosaguas, Pozuelo de Alarc\'on,
28223 Madrid
\\ \small Spain}


\maketitle

\begin{abstract}
In this paper a technique is given to recover the classical order of the method when explicit exponential Runge-Kutta methods integrate reaction-diffusion problems. Although methods of high stiff order for problems with vanishing boundary conditions can be constructed, that may imply increasing the number of stages and therefore, the computational cost seems bigger than the technique which is suggested here, which just adds some terms with information on the boundaries. Moreover, time-dependent boundary conditions are directly tackled here.
\end{abstract}


\msc{65M12 \and 65M20}

\section{Introduction}
\label{Introduction}

Order reduction is a well-known phenomenom that affects not only standard integrators with stages when integrating stiff problems but also exponential-type methods \cite{CR02,FOS,HO2}. The latter have the advantage of integrating the linear and stiff part of some differential problems in an exact way. This allows to construct methods which are linearly stable and \lq explicit' at the same time, which is not possible with standard integrators.

Different ideas to avoid order reduction with standard integrators and splitting methods (both exponential and not exponential) have been given in the literature but, up to our knowledge, order greater than 2 has not been achieved with them when integrating reaction-diffusion problems with time-dependent boundary conditions. (Look at \cite{acrnl,CBHW,EO,EO2,LV} for some of the results with splitting methods)

In this paper, we center on trying to increase accuracy among explicit exponential Runge-Kutta methods (EERK). Among them, Lawson ones are the easiest to construct since, given any standard Runge-Kutta (RK) method, there is directly a Lawson method associated to it through a change of variables given by an exponential function \cite{L}. However, some more methods were constructed in the literature considering other functions defined from exponentials and denoted by $\varphi_j$ ($j\ge 1$) \cite{HO2}. While Lawson methods show a severe order reduction when integrating initial boundary value problems \cite{acr1}, with some other EERK methods, the order reduction is not so big. Furthermore, in \cite{HO,LO0}, stiff order conditions are given when integrating semilinear parabolic problems with vanishing boundary conditions and then, some methods are constructed explicitly in order to achieve that order (see also \cite{LO}). However, imposing those conditions usually implies increasing the number of stages and thus, the computational cost.

For linear problems, a technique has already been suggested to avoid order reduction with both Lawson and other EERK methods \cite{acr2,CaM} not only in parabolic problems with vanishing boundary conditions, but also for time-dependent boundary values. This technique is based on applying the method of lines integrating firstly in time and then in space. Then,  in the first procedure, the exponential-type operators applied over functions with time-dependent boundary conditions are substituted by the solution of  suitable initial boundary value problems with appropriate boundary conditions which can be calculated exactly in terms of data of the problem.  The corresponding formulas after full discretization are also shown in \cite{acr2,CaM} and they are called exponential quadrature rules because, in fact, they do not depend at all on the stages of the method due to the fact that the problem is linear. The exponential RK quadrature rules in \cite{CaM} are in general more efficient than Lawson ones because, in the former, all the exponential-type operators are evaluated at the same argument and, because of that, their calculation through Krylov subroutines in \cite{niesen} is quicker \cite{CaM}. Moreover, in the same paper, for some particular methods, it is shown that the suggested technique to avoid order reduction with quadrature rules corresponding to few nodes is more efficient than considering the classical approach with higher stiff order because of a higher number of nodes.

For nonlinear problems and Lawson methods, a technique has also been suggested in \cite{CR}. The main difference with the linear case is that, in order to achieve local order $3$ and higher, numerical differentiation is required to approximate the necessary boundary values of the intermediate mentioned problems. Because of that, in order to prove convergence, a CFL condition of the form (\ref{cfl}) is required, although this condition is much less restrictive than that needed with standard explicit RK methods. For example, when $A$ is a second-order space derivative operator, $\gamma=2$ would be required in (\ref{cfl}) with the latter. However, with Lawson methods, $\gamma=1$ is just required with Dirichlet boundary conditions and no condition at all with Neumann ones because $\gamma=0$ in such a case. Moreover, the additional terms which need to be calculated in order to avoid order reduction seem to be very cheap to calculate with Krylov techniques \cite{CR}.

The aim of this paper is to see how to apply a similar technique with the EERK methods in \cite{HO}, for which $\varphi_j$ functions ($j\ge 1$) turn up, instead of just $\varphi_0(z)=e^{z}$, as it happens with Lawson methods. According to \cite{GG}, the more recent Krylov techniques are quicker with $\varphi_j$-operators when $j$ increases. Moreover, as distinct from Lawson methods, at least for the formula concerning the calculation of the numerical solution, the $\varphi_j$-operators ($j\ge0$) are always evaluated at the same argument. Furthermore, the same applies for the calculation of the stages when, in (\ref{a}), $r$ just takes the value $i$. As with exponential quadrature rules, this makes Krylov subroutines much quicker \cite{niesen} and, because of that, these methods are more efficient than Lawson ones. The analysis is much more complicated than with Lawson methods but the final conclusion is the same. Under a very mild CFL condition, order reduction can be avoided in a quite cheap way, since the additional calculations in fact just correspond to boundary values. That seems much cheaper than constructing methods with a high number of stages to increase the stiff order.

The paper is structured as follows. Firstly, a preliminaries section establishes the abstract framework for the problem to integrate and the required assumptions to deduce the main results in the theorems. Then, the appropriate time semidiscretization is suggested in Section 3, so that order reduction is avoided in the local error. In Section 4, the full discretization is performed considering some hypotheses on the space integrator. This includes simple finite differences, but also more complicated ones (which require a mass matrix, as the nine-point formula) and even finite elements. The formulas to be implemented are given in (\ref{Knih})-(\ref{Uhn}) and, afterwards, the result on the local error after full discretization is given, assuming that the required boundaries in (\ref{Knih})-(\ref{Uhn}) are calculated exactly. As this is not always possible in terms of data, Section 4 offers a first simplification of those boundaries in Tables 1 and 2, which guarantee that the optimal local error order after full discretization is conserved when searching for local order $2$, $3$ and $4$ in the timestepsize (respectively $p=1,2,3$ in Theorem \ref{thlocs}). Even with those simplifications, just some approximation of those boundaries can be given when $p=1$ with Neumann/Robin boundary conditions. Besides, when $p\ge 2$, with both Dirichlet and Neumann/Robin boundary conditions, numerical differentiation is required.  Considering then the error which comes from this approximation, the results on the final global error after full discretization is given in Subsection 4.6 for $p=1,2$ and $3$. Moreover, under the parabolic assumption (\ref{parabol}), it is justified that the global order on the timestepsize can be the same as the local one. Finally, in Section 5, some numerical experiments corroborate the previous results and show that the suggested technique to avoid order reduction with EERK methods of less stiff order can be more efficient than the classical approach with higher stiff order ones. As the proofs of the theorems are very technical, we have postponed the proofs to an appendix.

\section{Preliminaries}
\label{Preliminares}

Let $X$ and $Y$ be Banach spaces and let $A:D(A)\subset X \to X$ and
$\partial: X \to Y$ be linear operators. Our goal is to avoid order reduction when integrating in time through EERK methods the nonlinear abstract non homogeneous initial
boundary value problem
\begin{eqnarray}
\label{laibvp}
\begin{array}{rcl}
u'(t)&=&Au(t)+f(t,u(t)), \quad  0\le t \le T,\\
u(0)&=&u_0 \in X,\\
\partial u(t)&=&g(t)\in Y, \quad  0\le t \le T.
\end{array}
\end{eqnarray}

The abstract setting (\ref{laibvp}) permits to cover a wide range
of nonlinear evolutionary problems governed by partial
differential equations. We use the following hypotheses, similar to the ones in \cite{acrnl} when avoiding the same kind of problems with exponential splitting methods.

\begin{enumerate}
\item[(A1)] The boundary operator $\partial:D(A)\subset X\to Y$ is
onto and $g\in C^1([0,T],Y)$.

\item[(A2)]  Denoting $A_0:D(A_0)=Ker(\partial)\subset X \to X$, the restriction of $A$
to Ker($\partial$), at least one of these assumptions is satisfied:
\begin{enumerate}
\item
Ker($\partial$) is dense in $X$ and $A_0$ is the infinitesimal generator of a $C_0$-
semigroup $\{e^{t A_0}\}_{t\ge 0}$ in $X$ of negative type $\omega$.
\item
$D(A)$ is dense in $X$ and $A_0$ generates a bounded holomorphic semigroup $\{e^{t A_0}\}_{t\ge 0}$ in $X$ of negative type $\omega$.
\end{enumerate}
\item[(A3)] If $z \in \mathbb{C}$ satisfies $\Re (z) >\omega$ and $v
\in Y$, then the steady state problem
\begin{eqnarray}
 Ax &=& zx,  \nonumber \\
 \partial x&=&v,
\nonumber
\end{eqnarray}
possesses a unique solution denoted by $x=K(z)v$. Moreover, the
linear operator $K(z): Y \to D(A)$ satisfies
\begin{eqnarray}
\label{stationaryoperator} \| K(z)v\| \le C\|v\|,
\end{eqnarray}
where the constant $C$ holds for any $z$ such that $Re (z) \ge
\omega_0 > \omega$.

\item[(A4)] The nonlinear source $f$ belongs to $C^1([0,T] \times
X, X)$.

\item[(A5)] The solution $u$ of (\ref{laibvp}) satisfies $u\in
C^1([0,T], X)$, $u(t) \in D(A)$ for all $t \in [0,T]$ and $Au \in C([0,T], X)$.

\end{enumerate}

The well-posedness of problem (\ref{laibvp}) is assured by (A1)-(A4), as explained in
\cite{acrnl}. Besides, as it was also justified there, (A4) is quite restrictive if $X=L^p(\Omega)$  with $\Omega$ a bounded domain in $\mathbb{R}^d$. However, if the supremum norm is chosen, (A4) is satisfied whenever
$f$ has the form
\begin{eqnarray}
\label{nonlinearterm} f(t,u)= \phi (u) + h(t),
\end{eqnarray}
with $\phi\in C^1(\mathbb{C},\mathbb{C})$ and  $h\in C^1([0,T], X)$.
For simplicity, we will assume from now on that $f$ has the form (\ref{nonlinearterm}), apart from the following hypotheses:
\begin{enumerate}
\item[(A6)] There exists a natural value $\bar{m}(A)\ge 1$ such that,  whenever $w\in D(A^l)$  and $\phi \in C^{m+ \bar{m}(A)l}(\mathbb{C},\mathbb{C})$ for natural $l$ and $m$, $\phi^{(m)}(w)\in D(A^l)$.

\item[(A7)] For every natural $l\ge 0$, there exists a norm $\|\cdot \|_l$ in $D(A^l)\subset X$ such that, for every natural $m\ge 1$ and $u,v\in D(A^l)$, whenever $\phi \in C^{m+ \bar{m}(A)l}(\mathbb{C},\mathbb{C})$, $\phi^{(m)}(u) v^m \in D(A^l)$ and
    $$
    \|A^l \big[\phi^{(m)}(u) v^m]\big]\| \le C(\phi,u) \|v\|_l^m,
    $$
    for some constant $C(\phi,u)$ which depends on $\phi^{(m)}(u),\dots,\phi^{(m+ \bar{m}(A)l)}(u)$ and $\|u\|_l$.
\item[(A8)] For every natural $l\ge 1$ and the norm $\|\cdot \|_l$ in (A7), for every $u,v\in D(A^l)$, whenever $\phi \in C^{(1+ \bar{m}(A)l)}(\mathbb{C},\mathbb{C})$,
    $$
    \|A_0^{-1} A^l \big[\phi'(u)v\big]\| \le \bar{C}(\phi,u) \|v\|_{l-1},
    $$
    for some constant $\bar{C}(\phi,u)$ which depends on $\phi'(u),\dots,\phi^{(1+ \bar{m}(A)l)}(u)$ and $\|u\|_l$.

\end{enumerate}
Notice that, when $A$ is a space differential operator, (A6) means that $\phi$ and some of its derivatives do not reduce the regularity in space of their arguments for regular enough $\phi$. As for (A7), it assures that Taylor expansions with the unbounded operator $A$ can be performed when there is enough regularity. Finally, (A8) is used to apply a summation-by-parts argument which explains that the local order is the same as the global one. (These assumptions will be justified for the example in the numerical experiments.)

In the remaining of the paper, we always suppose that (A1)-(A8)
are satisfied. However,  more
regularity will be assumed in certain results.

Because of hypothesis (A2), $\{\varphi_j(t A_0)\}_{j=0}^{\infty}$   are
bounded operators for $t>0$, where $\{\varphi_j\}$ are the
standard functions which are used in exponential methods \cite{HO2} and which
are defined by
\begin{eqnarray}
\varphi_j( t A_0)=\frac{1}{t^j} \int_0^t
e^{(t-\tau)A_0}\frac{\tau^{j-1}}{(j-1)!}d\tau, \quad j \ge 1.
\label{varphi}
\end{eqnarray}
It is well-known that they can be calculated in a recursive way through the formulas
\begin{eqnarray}
\varphi_{j+1}(z)=\frac{\varphi_j(z)-1/j!}{z}, \quad z \neq 0,
\qquad \varphi_{j+1}(0)=\frac{1}{(j+1)!}, \qquad \varphi_0(z)=e^z.
\label{recurf}
\end{eqnarray}

For the time integration, we will center on EERK methods which, when applied to a
finite-dimensional nonlinear problem  like
\begin{eqnarray}
U'(t) = M U(t)+F(t,U(t)), \label{linfd}
\end{eqnarray}
with $M$ a certain matrix, read like this at each step
\begin{eqnarray}
K_{n,i}&=&e^{c_i k M}U_n+k \sum_{j=1}^{i-1} a_{ij}(k M) F(t_n+c_j k, K_{n,j}), \quad i=1,\dots,s, \label{etapas} \\
U_{n+1}&=&e^{k M}U_n+k \sum_{i=1}^s b_i(k M)F(t_n+c_i k,K_{n,i}).
\label{eerk}
\end{eqnarray}
Here,  $k >0$ is the time stepsize, $t_n=t_0+nk$ and, in principle, $a_{ij}(k M)$ and $b_i(k M)$ are linear combinations of $\{ \varphi_l(c k M)\}_{l=1}^s$, for some real value $c$. More precisely, we will assume that
\begin{eqnarray}
a_{ij}(z)&=&\sum_{l,r=1}^s \lambda_{i,j,l,r} \varphi_l( c_r z), \label{a} \\
b_i(z)&=& \sum_{l=1}^s \mu_{i,l} \varphi_l(z), \nonumber
\end{eqnarray}
for some constants $\lambda_{i,j,l,r}$ and $\mu_{i,l}$. We notice that, usually, in the formula for $a_{ij}$, $l$ just varies from $1$ to $i-1$ and $r$ just takes the value $i$, but we consider the more general case because it does not imply any difficulty in the analysis. We also notice that, usually, the underlying method satisfies $\sum_{j=1}^s a_{ij}(0)=c_i$, which implies that $\sum_{j,l,r=1}^s \lambda_{i,j,l,r}/l!=c_i$. We will assume this from now on for the simplicity of some expressions in Section \ref{SRB}.

\section{Suggestion for the time semidiscretization}

As already mentioned in the introduction, we will integrate (\ref{laibvp}) firstly in time in order to avoid order reduction. However,  (\ref{laibvp}) is not an ODE system like (\ref{linfd}), but an initial boundary value problem. Because of that,
for each of the terms in (\ref{etapas})-(\ref{eerk}), we must suggest the solution of a suitable initial boundary value problem. Simulating what was already done with exponential functions and Lawson methods \cite{acr2},  whenever $\alpha \in X$ and $\beta\in D(A^{r+1})$, we consider the generalised solution of this problem, which we will denote by $\phi_{0,r,\alpha,\beta}(\tau)$,
\begin{eqnarray}
v'(\tau)&=& A v(\tau), \nonumber \\
v(0)&=& \alpha, \nonumber \\
\partial v (\tau) &=& \partial (\sum_{l=0}^r \frac{\tau^l}{l!} A^l \beta). \label{phi0}
\end{eqnarray}
On the other hand, in a similar way to what was done for $\varphi_j$-functions and exponential quadrature rules \cite{CaM}, and under the same assumptions on $\alpha$ and $\beta$, we also consider the solution of this other problem, which we will denote by $\phi_{j,r,\alpha,\beta}(\tau)$,
\begin{eqnarray}
v'(\tau)&=& (A-\frac{j}{\tau})v(\tau)+\frac{1}{(j-1)!\tau}\alpha, \nonumber \\
v(0)&=& \frac{1}{j!}\alpha, \nonumber \\
\partial v (\tau) &=& \partial (\sum_{l=0}^r \frac{\tau^l}{(l+j)!} A^l \beta). \label{phij}
\end{eqnarray}
As it was already proved in Lemma 3.1 of \cite{acr2} for $\phi_{0,r,\alpha,\beta}$, and as it can be proved in a similar way for $\phi_{j,r,\alpha,\beta}$, considering Lemmas 6 and 7 in \cite{CaM}, for $r\ge 0$,
\begin{eqnarray}
\phi_{0,r,\alpha,\beta}(\tau)&=&\sum_{l=0}^r \frac{\tau^l}{l!} A^l \beta+e^{\tau A_0}(\alpha-\beta)+\tau^{r+1} \varphi_{r+1}(\tau A_0) A^{r+1} \beta, \label{phi0r} \\
\phi_{j,r,\alpha,\beta}(\tau)&=&\sum_{l=0}^r \frac{\tau^l}{(l+j)!} A^l \beta+\varphi_j(\tau A_0)(\alpha-\beta)+\tau^{r+1} \varphi_{j+r+1}(\tau A_0) A^{r+1} \beta. \label{phijr}
\end{eqnarray}
Then,  when the EERK method has non-stiff order $p$, our suggestion in order to avoid order reduction, when the problem is regular enough, is to consider
as stages, from the previous numerical value $u_n$,
\begin{eqnarray}
K_{n,i}= \phi_{0,p-1,u_n,u(t_n)}(c_i k)+k \sum_{j=1}^{i-1} \sum_{l,r=1}^s \lambda_{i,j,l,r} \phi_{l,p-2,F_{n,j}, \bar{F}_{n,j}}(c_r k), \label{etsro}
\end{eqnarray}
where
\begin{eqnarray}
F_{n,j}=f(t_n+c_j k, K_{n,j}), \quad \bar{F}_{n,j}=f(t_n+c_j k, \bar{K}_{n,j}), \label{Fn}
\end{eqnarray}
 with $\bar{K}_{n,j}$  obtained as $K_{n,j}$ but substituting $u_n$ by $u(t_n)$ and, recursively, $F_{n,l}$ by $\bar{F}_{n,l}$. Then, we suggest as the numerical solution after timestepsize $k$,
\begin{eqnarray}
u_{n+1}=\phi_{0,p,u_n,u(t_n)}(k)+k \sum_{i=1}^s \sum_{l=1}^s \mu_{i,l} \phi_{l,p-1,F_{n,i},\bar{F}_{n,i}}(k). \label{usro}
\end{eqnarray}
We notice that, when $p=1$, instead of (\ref{etsro}), it should be used
\begin{eqnarray}
K_{n,i}=\phi_{0,0,u_n,u(t_n)}(c_i k)+k \sum_{j=1}^{i-1} \sum_{l,r=1}^s \lambda_{i,j,l,r} \varphi_l(c_r k A_0)F_{n,j},
\label{Kni1}
\end{eqnarray}
since $\phi_{l,p-2,F_{n,j}, \bar{F}_{n,j}}$ should be the solution of (\ref{phij}), with vanishing boundary conditions.

\subsection{Local error}
\label{localerror}

Before stating and proving the result on the local error which assures that order reduction would be avoided in such a way under enough assumptions of regularity, we remind that the conditions for non-stiff order $p$ come from considering the Taylor series of $U(t_{n+1})$ in terms of $U(t_n)$ when $U(t)$ is the solution of the ordinary differential system (\ref{linfd}), i.e.,
\begin{eqnarray}
U(t_{n+1})&=& U(t_n)+ k[M U(t_n)+F_n]\nonumber \\
&&+\frac{k^2}{2}[M^2 U(t_n)+M F_n +F_{t,n}+F_{u,n} M U(t_n)+F_{u,n} F_n] \nonumber \\
&&+\sum_{j=3}^p\frac{k^j}{j!} \sum_{m=1}^{n(j)} \gamma_{m,j} ED_{m,j,n} (F,M)(t_n,U(t_n))+O(k^{p+1}), \label{taylorUord}
\end{eqnarray}
where $F_{*,n}$ denotes the corresponding evaluation of $F$ or some of its derivatives when evaluated at $(t_n,U(t_n))$ and $\{ED_{m,j,n}(F,M)\}_{m=1}^{n(j)}$ are the independent elementary differentials associated to the power $k^j$ ($0\le j\le p$) and also evaluated at $(t_n,U(t_n))$ \cite{HNW}.  The order conditions come from equating those coefficients $\gamma_{m,j}$ with the ones corresponding to the Taylor series of the corresponding numerical solution (\ref{etapas})-(\ref{eerk}) starting from $U(t_n)$. More precisely, the result of developing in Taylor series around $(t_n,U(t_n))$,
\begin{eqnarray}
\bar{U}_{n+1}&=& \sum_{l=0}^p \frac{k^l}{l!} M^l U(t_n) \nonumber \\
&&+k \sum_{i=1}^s \sum_{l=1}^s \mu_{i,l} \sum_{j=0}^{p-1} \frac{k^j}{(j+l)!} M^j F(t_n+c_i k, \bar{K}_{n,i})+O(k^{p+1}),
\label{Ubord}
\end{eqnarray}
where, recursively,
\begin{eqnarray}
\bar{K}_{n,i}&=&\sum_{l=0}^{p-1} \frac{(c_i k)^l}{l!} M^l U(t_n) \nonumber \\
&&+k\sum_{j=1}^{i-1} \sum_{l,r=1}^s \lambda_{i,j,l,r} \sum_{ll=0}^{p-2} \frac{(c_r k)^{ll}}{(ll+l)!} M^{ll} F(t_n+c_j k,\bar{K}_{n,j})+O(k^p). \label{Kbord}
\end{eqnarray}

For the numerical integration of (\ref{laibvp}) through (\ref{etsro})-(\ref{usro}), we then have the following result for the local error $\rho_n=\bar{u}_{n+1}-u(t_{n+1})$, where $\bar{u}_{n+1}$ is defined also through (\ref{etsro})-(\ref{usro}), but starting from $u(t_n)$ instead of $u_n$.

\begin{theorem}
Under hypotheses (A1)-(A8), if the explicit exponential Runge-Kutta method (\ref{etapas})-(\ref{eerk}) has non-stiff order $\ge p$, when integrating (\ref{laibvp}) through (\ref{etsro})-(\ref{usro}) with
$u\in C([0,T], D(A^{p+1}))\cap C^{p+1}([0,T],X)$, $\phi\in C^{1+(p-1)\bar{m}(A)}(\mathbb{C},\mathbb{C})$,  $h\in C^{p}([0,T],D(A^p)),$
it happens that the local error satisfies $\rho_n=O(k^{p+1})$.

Moreover, if in formulas (\ref{etsro})-(\ref{usro}), we had substituted $p$ by $\hat{p}$ with $1\le \hat{p}\le p-1$, whenever $u\in C([0,T], D(A^{\hat{p}+2}))\cap C^{\hat{p}+2}([0,T],X)$, $\phi\in C^{1+\hat{p}\bar{m}(A)}(\mathbb{C},\mathbb{C})$ and $h\in C^{\hat{p}+1}([0,T],D(A^{\hat{p}+1}))$, it happens that not only the local error satisfies $\rho_n=O(k^{\hat{p}+1})$ but also $A_0^{-1} \rho_n=O(k^{\hat{p}+2})$.
\label{thloc}
\end{theorem}

The result of the previous theorem corresponding to $\hat{p}< p$ is interesting because, for the non-stiff order $p$, some times the corresponding boundaries in (\ref{phi0}) and (\ref{phij}) cannot be calculated exactly in terms of the data of the problem. We can always resort to numerical differentiation to approximate them, but the more terms which need to be calculated in such a way, the more complicated and unstable the algorithm may become. That will be studied in detail afterwards.

\section{Suggestion for the full discretization}
\label{suggestionfulldiscretization}

As distinct from the standard approach when integrating partial differential problems, we suggest to discretize firstly in time as described in the previous section and then to integrate in space the resulting auxiliary problems. For the suggestion of a space discretization, we must first fix a certain space $X$. As in \cite{acrnl} and \cite{CR}, we take
$X=C(\overline{\Omega})$ for a certain bounded domain $\Omega\in \mathbb{R}^d$ and the maximum norm.

\subsection{Space discretization}

In this subsection we consider  a
certain grid $\Omega_h$ (of $\Omega$) over which the approximated
numerical solution will be defined. Therefore, the numerical
approximation belongs to $\mathbb{C}^N$, where $N$ is the number of nodes
in the grid, and the maximum norm $\|\cdot\|_h$ in this space will also be considered.
Notice that, with Dirichlet boundary
conditions, nodes on the boundary are not part of the approximation while, when
using Neumann or Robin boundary conditions, the approximation at those nodes must be considered.

We will make similar hypotheses for the space discretization as those in \cite{acrnl} and  \cite{CaM}. We state them here for the sake of clarity.
We denote the projection operator
\begin{eqnarray}
 P_h : X \to \mathbb{C}^N, \nonumber
\end{eqnarray}
to that which takes a function to its values over the grid $\Omega_h$. Notice that the following bound then holds
\begin{eqnarray}
\|P_h u \|_{h} \le \|u\|_{h}. \label{acotp}
\end{eqnarray}
 Then, we assume that the discretization of
 $$A u =F, \qquad \partial u=g,$$
 is given by
\begin{eqnarray}
A_{h,0}U_h+C_h g=P_h F+ D_h \partial F,
\label{spacediscr}
\end{eqnarray}
where $A_{h,0}$ is the matrix which discretizes $A_0$ and $C_h, D_h: Y
\to \mathbb{C}^N$ are other operators associated to the discretization of $A$, which take into account the information on the boundary of $u$ and $F$.

Notice that, as $f$ is given by (\ref{nonlinearterm}), it has sense as a
function from $[0,T]\times \mathbb{C}^N$ on $\mathbb{C}^N$ and,
for each $t \in [0,T]$ and $u \in X$,
\begin{eqnarray}
\label{fandph} P_hf(t,u)= f(t, P_hu).
\end{eqnarray}

Moreover, we consider the following hypotheses:
\begin{enumerate}
\item[(H1)] The matrix $A_{h,0}$ satisfies
\begin{enumerate}
\item $A_{h,0}$ is invertible and $\|A_{h,0}^{-1}\|_h \le C$ for
some $h$-independent constant $C$.
\item $\|\varphi_j(tA_{h,0})\|_h \le C_j, \quad j=0,\dots,s, \quad t\in [0,T]$.
\end{enumerate}
\item[(H2)] We define the elliptic projection $R_{h}:D(A) \to
\mathbb{C}^N$  as the solution of
\begin{eqnarray}
A_{h,0} R_{h} u+ C_{h}\partial u=P_h Au +D_h \partial A u. \label{rh}
\end{eqnarray}
\begin{enumerate}
\item[(a)] There exists a subspace $Z \subset D(A)$ such
that, for $u \in Z$,
$A_0^{-1} u \in Z$ and $e^{tA_0}u\in Z$, for $t \in [0,T]$. Besides,
there exists a natural value $\bar{m}(Z)\ge 1$ such that, whenever $w\in D(A^l)$ for natural $l$ and $A^l w\in Z$, it happens that $A^l \phi(w)\in Z$
if $\phi \in C^{\bar{m}(Z)+l\bar{m}(A)}(\mathbb{C},\mathbb{C})$. Moreover,
for some $\varepsilon_{h}$ and $\eta_{h}$ which are
both small with $h$, if $u\in Z$,
\begin{equation}
\label{consistency}
 \left\| A_{h,0}({P_hu-R_{h}u})
\right\|_h \le \varepsilon_{h} \left\| u \right\|_Z, \quad
\left\|P_hu-R_{h}u \right\|_h \le \eta_{h} \left\| u \right\|_Z.
\end{equation}
(Although obviously, because of (H1a), $\eta_h$ could be taken as
$C\varepsilon_h$, for some discretizations $\eta_h$ can decrease more quickly with $h$
than $\varepsilon_h$ and that leads to better error bounds in the
following sections.)
\item[(b)] $\|A_{h,0}^{-1} C_h\|_h \le C''$ for some constant
$C''$ which does not depend on $h$. This resembles the continuous
maximum principle which is satisfied because of
(\ref{stationaryoperator}) when $z=0$.
\item[(c)] $\|D_h\|_h$ is uniformly bounded on $h$.
\end{enumerate}
\item[(H3)] The nonlinear source $f$ belongs to
$C^1([0,T] \times \mathbb{C}^N, \mathbb{C}^N)$ and the derivative with respect to the variable in $\mathbb{C}^N$
is uniformly bounded in a neighbourhood of the solution where the numerical approximation stays.
\end{enumerate}

\subsection{Full discretization formulas depending on $\partial A^l u(t_n)$ and $\partial A^l \bar{F}_{n,i}$}

As in (\ref{etsro}) and (\ref{usro}), $\alpha=u_n$ in the part which corresponds to (\ref{phi0}) and $\alpha=F_{n,i}$ in that corresponding to (\ref{phij}), we suggest to discretize (\ref{phi0}) and (\ref{phij}) in space considering (\ref{spacediscr}) and taking as initial condition for the space discretized system
$U_h^n$ for (\ref{phi0}) with boundary $\partial u(t_n)$ and $F_{n,i,h}$ for (\ref{phij}) with boundary $\partial \bar{F}_{n,i}$.

More precisely, for $r=p-1$ or $r=p$, depending on whether we are calculating the stages or the solution, we consider the space discretized system related to (\ref{phi0}):
\begin{eqnarray}
V_h'(\tau)&=&A_{h,0} V_h(\tau)+C_h \partial[\sum_{l=0}^r \frac{\tau^l}{l!} A^l u(t_n)]-D_h \partial[\sum_{l=1}^r \frac{\tau^{l-1}}{(l-1)!} A^l u(t_n)], \nonumber \\
V_h(0)&=&U_h^n. \nonumber
\end{eqnarray}
By the variation-of-constants formula and the definition of $\varphi_l$ (\ref{varphi}), it is clear that
\begin{eqnarray}
\phi_{0,r,U_h^n,u(t_n)}^h(\tau)&\equiv&e^{\tau A_{h,0}}U_h^n +\sum_{l=0}^{r-1} \tau^{l+1} \varphi_{l+1} (\tau A_{h,0})[C_h \partial A^l u(t_n)-D_h \partial A^{l+1}u(t_n)]\nonumber \\
&&+\tau^{r+1} \varphi_{r+1}(\tau A_{h,0}) C_h \partial A^r u(t_n)
\nonumber
\end{eqnarray}
solves this problem. (We consider this notation for similarity with the continuous one in (\ref{phi0r}).)
In a similar way, for the space discretized system related to (\ref{phij}), taking $r=p-2$ or $r=p-1$ depending on whether we are calculating the stages or the solution, we consider
\begin{eqnarray}
V_h'(\tau)&=&(A_{h,0}-\frac{j}{\tau}) V_h(\tau)+\frac{1}{(j-1)! \tau} F_{n,i,h} +C_h \partial[\sum_{l=0}^r \frac{\tau^l}{(l+j)!} A^l \bar{F}_{n,i}] \nonumber \\
&& -D_h \partial [\sum_{l=1}^r \frac{l \tau^{l-1}}{(l+j)!} A^l \bar{F}_{n,i}]-j D_h \partial [\sum_{l=0}^r \frac{\tau^{l-1}}{(l+j)!} A^l \bar{F}_{n,i}]+\frac{1}{(j-1)! \tau} D_h \partial \bar{F}_{n,i}, \nonumber
\end{eqnarray}
which can be written like this for $r>-1$,
\begin{eqnarray}
V_h'(\tau)&=&(A_{h,0}-\frac{j}{\tau}) V_h(\tau)+\frac{1}{(j-1)! \tau} F_{n,i,h} \nonumber \\
&&+\sum_{l=0}^{r-1} \frac{\tau^l}{(l+j)!}[C_h \partial A^l \bar{F}_{n,i}-D_h \partial A^{l+1} \bar{F}_{n,i}]+\frac{\tau^r}{(r+j)!}C_h \partial A^r \bar{F}_{n,i} \nonumber \\
V_h(0)&=&\frac{1}{j!}F_{n,i,h}. \nonumber
\end{eqnarray}
Then, the solution of this, in the same way that it was argued for exponential quadrature rules in \cite{CaM}, can be written like this for $\tau>0$ and $r>-1$,
\begin{eqnarray}
\phi_{j,r,F_{n,i,h},\bar{F}_{n,i}}^h(\tau)&\equiv&\varphi_j (\tau A_{h,0}) F_{n,i,h}+\sum_{l=0}^{r-1} \tau^{l+1} \varphi_{j+l+1} (\tau A_{h,0}) [C_h \partial A^l \bar{F}_{n,i}- D_h \partial A^{l+1} \bar{F}_{n,i}]
\nonumber\\
&& +\tau^{r+1} \varphi_{j+r+1}(\tau A_{h,0})C_h \partial A^r \bar{F}_{n,i}. \nonumber
\end{eqnarray}
Considering all this, the suggestion for the full discretization formulas is, in principle:
\begin{eqnarray}
\lefteqn{K_{n,i,h}=e^{c_i k A_{h,0}} U_h^n+\sum_{l=0}^{p-2} (c_i k)^{l+1} \varphi_{l+1}(c_i k A_{h,0}) [C_h \partial A^l u(t_n)-D_h \partial A^{l+1} u(t_n)]} \nonumber \\
&&\hspace{3cm}+(c_i k)^p \varphi_p(c_i k A_{h,0})C_h \partial A^{p-1} u(t_n) \nonumber \\
&&+k \sum_{j=1}^{i-1} \sum_{l,r=1}^s \lambda_{i,j,l,r}[\varphi_l(c_r k A_{h,0}) F_{n,j,h}+\sum_{ll=0}^{p-3}(c_r k)^{ll+1} \varphi_{l+ll+1}(c_r k A_{h,0})[C_h \partial A^{ll} \bar{F}_{n,j}-D_h \partial A^{ll+1} \bar{F}_{n,j}] \nonumber \\
&&\hspace{3cm}+(c_r k)^{p-1} \varphi_{l+p-1}(c_r k A_{h,0})C_h \partial A^{p-2} \bar{F}_{n,j}],
 \label{Knih} \\
\lefteqn{U_h^{n+1}=e^{k A_{h,0}}U_h^n+\sum_{l=0}^{p-1} k^{l+1} \varphi_{l+1}(k A_{h,0})[C_h \partial A^l u(t_n)-D_h \partial A^{l+1} u(t_n)]} \nonumber \\
&&\hspace{3cm}+k^{p+1} \varphi_{p+1} (k A_{h,0})C_h \partial A^p u(t_n) \nonumber \\
&&+k \sum_{i=1}^s \sum_{l=1}^s \mu_{i,l}\bigg[\varphi_l(k A_{h,0})F_{n,i,h}+\sum_{ll=0}^{p-2} k^{ll+1} \varphi_{l+ll+1}(k A_{h,0})[C_h \partial A^{ll} \bar{F}_{n,i}- D_h \partial A^{ll+1} \bar{F}_{n,i}] \nonumber \\
&&\hspace{3cm}+k^p \varphi_{l+p}(k A_{h,0}) C_h \partial A^{p-1} \bar{F}_{n,i}\bigg], \label{Uhn}
\end{eqnarray}
where $F_{n,i,h}=f(t_n+c_i k, K_{n,i,h})$ and where we notice that, for $p=1$, due to (\ref{Kni1}), (\ref{Knih}) should read
\begin{eqnarray}
K_{n,i,h}&=&e^{c_i k A_{h,0}} U_h^n+c_i k \varphi_{1}(c_i k A_{h,0}) C_h \partial  u(t_n)\nonumber \\
&&+k \sum_{j=1}^{i-1} \sum_{l,r=1}^s \lambda_{i,j,l,r} \varphi_l(c_r k A_{h,0}) F_{n,j,h}. \label{Knih1}
\end{eqnarray}

\subsection{Local error of the full discretization}

We define the local error of the full discretization by $\rho_{n,h}=\bar{U}_h^{n+1}-P_h u(t_{n+1})$, where $\bar{U}_h^{n+1}$ is defined through (\ref{Knih})-(\ref{Uhn}) but substituting $U_h^n$ by $P_h u(t_n)$. The following theorem states how this error behaves.

\begin{theorem} Under the first set of hypotheses of Theorem \ref{thloc}, (H1)-(H3) and assuming also that $\phi \in C^{\bar{m}(Z)+p \bar{m}(A)}(\mathbb{C},\mathbb{C})$,
\begin{eqnarray}
A^l u \in C([0,T], Z), \quad l=0,\dots,p+1, \quad A^l h\in C([0,T],Z),  \quad l=0,\dots,p,
\label{regfle}
\end{eqnarray}
it happens that $\rho_{n,h}=O(k^{p+1}+k \varepsilon_h)$.

Moreover, under the second set of hypotheses of Theorem \ref{thloc},(H1)-(H3) and assuming also that $\phi \in C^{\bar{m}(Z)+\hat{p} \bar{m}(A)}(\mathbb{C},\mathbb{C})$ and (\ref{regfle}), but with $\hat{p}$ instead of $p$, it happens that, not only $\rho_{n,h}=O(k^{\hat{p}+1}+k \varepsilon_h)$ but also $A_{h,0}^{-1} \rho_{n,h}=O(k^{\hat{p}+2}+k \eta_h+k^2\varepsilon_h)$.
\label{thlocfd}
\end{theorem}

\subsection{Simplification of the required boundaries without losing order}
\label{SRB}

As it also happened with Lawson methods and nonlinear problems \cite{CR}, the exact calculation of the required boundaries is not usually possible in terms of the data of the problem and just some numerical approximation of those values can be considered. Because of that, before considering any numerical approximation of them, we firstly suggest some simplifications of the previously suggested boundary values which do not alter the local error after full discretization. More precisely, for $p=1,2,3$ in (\ref{Knih})-(\ref{Uhn}), we suggest some simplifications which do not alter the respective order $p+1$ for the local error after full discretization, which has just been proved in the previous subsection. More precisely, we suggest the simplifications for $\partial A^l \bar{F}_{n,i}$ in Tables \ref{t1} and \ref{t2}, which come from developing into Taylor series $\bar{K}_{n,i}$ in (\ref{Kbord}) around $t_n$. (Notice that the simplification is different depending on whether the term is in the stages (\ref{Knih}) or in the approximation (\ref{Uhn})).

In the following, we will state more precisely under which assumptions the order of the  full discretization local error after that simplification (which we will denote by $\rho_{n,h}^{simp}$) will be conserved.

\begin{table}[t]
\caption{Simplifications for the boundaries in (\ref{Knih})}
\label{t1}
\begin{tabular}{cc}
\hline\noalign{\smallskip}
$p$ & Simplification in $K_{n,i,h}$  \\ \noalign{\smallskip}\hline\noalign{\smallskip}
2 & $\partial \bar{F}_{n,i} \approx \partial f(t_n,u(t_n))$   \\ \hline\noalign{\smallskip}
3 & $\begin{array}{rcl} \partial \bar{F}_{n,i} &\approx& \partial f (t_n+c_i k,u(t_n) +c_i k  \dot{u}(t_n) ) \\
\partial A \bar{F}_{n,i} &\approx& \partial A f(t_n,u(t_n))
\end{array}$ \\ \hline\noalign{\smallskip}
\end{tabular}
\end{table}

\begin{table}[t]
\caption{Simplifications for the boundaries in (\ref{Uhn})}
\label{t2}
\begin{tabular}{cc}
\hline\noalign{\smallskip} $p$
& Simplification in $U_h^{n+1}$
\\ \noalign{\smallskip}\hline\noalign{\smallskip}
1 &  $\partial \bar{F}_{n,i} \approx \partial f(t_n,u(t_n))$ \\ \noalign{\smallskip}\hline\noalign{\smallskip}
2 & $\begin{array}{rcl} \partial \bar{F}_{n,i} &\approx& \partial f (t_n+c_i k,u(t_n) +c_i k \dot{u}(t_n)) \\
\partial A \bar{F}_{n,i} &\approx& \partial A f(t_n,u(t_n))
\end{array}$ \\ \noalign{\smallskip}\hline\noalign{\smallskip}
3 & $\begin{array}{rcl} \partial \bar{F}_{n,i} &\approx& \partial f \bigg(t_n+c_i k,u(t_n) +k c_i A u(t_n)+\frac{(c_i k)^2}{2} A^2 u(t_n) \\
& & \hspace{1cm}+k\sum_{j,l,r} \lambda_{i,j,l,r}\big[ \frac{1}{l!}f \big(t_n+c_j k,u(t_n) +c_j k \dot{u}(t_n))+\frac{c_r k}{(l+1)!} A f(t_n,u(t_n))\big]\bigg) \\
\partial A \bar{F}_{n,i} &\approx& \partial A f (t_n+c_i k,u(t_n) +c_i k \dot{u}(t_n)) \\
\partial A^2 \bar{F}_{n,i} &\approx& \partial A^2 f(t_n,u(t_n))
\end{array}$ \\ \hline\noalign{\smallskip}
\end{tabular}

\end{table}

\begin{theorem}
Whenever $p=1,2,3$, under hypotheses (A1)-(A8) and (H1)-(H3), if the EERK method (\ref{etapas})-(\ref{eerk}) has non-stiff order $\ge p$, when integrating (\ref{laibvp}) through (\ref{Knih})-(\ref{Uhn}) with $u\in C([0,T], D(A^{p+1}))\cap C^{p+1}([0,T],X)$, $\phi\in C^{\bar{m}(Z)+p\bar{m}(A)}(\mathbb{C},\mathbb{C})$, $h\in C^p([0,T],D(A^p))$ and (\ref{regfle}), it happens that the local error of the full discretization with the simplified boundaries for $U_h^{n+1}$ in Tables \ref{t1} and \ref{t2} satisfies $\rho_{n,h}^{simp}=O(k^{p+1}+k\varepsilon_h)$.

Furthermore, if the method has non-stiff order $\ge p+1$,  $u\in C([0,T], D(A^{p+2}))\cap C^{p+2}([0,T],X)$ and $h\in C^{p+1}([0,T],D(A^{p+1}))$
and, just for the case   $p\geq 2$, the following bound holds
\begin{eqnarray}
\|A_0^{-1} [\phi'(u(t))A_0 w] \|\le C \|w\|, \mbox{ for every }t\in [0,T], \, w\in D(A_0), \label{cspc}
\end{eqnarray}
it happens that $A_{h,0}^{-1} \rho_{n,h}^{simp}=O(k^{p+2}+k \eta_h+k^2 \varepsilon_h)$.
\label{thlocs}
\end{theorem}

\subsection{Calculation of the required boundaries}
In this subsection, we will state the error which is committed when approximating not only the corresponding simplified boundaries in (\ref{Knih})-(\ref{Uhn}) described in Tables \ref{t1} and \ref{t2}, but also the error in the approximation  of $\partial A^l u(t_n)$ for the different values of $l$. The latter and part of the former was already studied in detail in \cite{CR} for Dirichlet and Robin/Neumann boundary conditions when considering Lawson methods. There are some differences in the arguments where the terms of the form $\partial A^l f$ must be evaluated now, but the error which is committed when approximating these boundaries in terms of the data of the problem (\ref{laibvp}), through the numerical approximation of the problem itself, or using numerical differentiation, is very similar. We state it here without entering into the details but reminding that $e_{n,h}=U_h^n -P_h u(t_n)$ (and therefore a bound for this error turns up when substituting $u(t_n)$ at some node in the boundary by the numerical approximation itself), $\mu_{k,1}$ is the norm of the error of some numerical differentiation formula for approximating the first derivative in time, $\nu_h$ the norm of the error of some numerical differentiation formula for approximating a certain $\gamma$th-derivative in space and $\mu_{k,2}$ the norm of the error of some numerical differentiation formula for approximating the second derivative in time.

For simplicity, we will separate the cases $p=1,2,3$ (which correspond respectively to Tables \ref{t3},\ref{t4} and \ref{t5}), distinguishing between the error which is committed for every term in the boundaries for the stages and for the numerical solution, and taking into account the term to be approximated in formulas (\ref{Knih})-(\ref{Uhn}) through their simplifications in Tables \ref{t1} and \ref{t2}. We also notice that, in 1-dimensional problems, when $A=D_{xx}$, as justified in \cite{CR} for $p=2,3$, $\gamma=1$ with Dirichlet boundary conditions and $\gamma=0$ for Robin/Neumann ones.

\begin{table}
\caption{Errors which are committed at each step when approximating the corresponding simplifications of the boundary terms in Tables \ref{t1} and \ref{t2}  to achieve local order $2$}
\label{t3}
\begin{tabular}{llcc}
\hline\noalign{\smallskip}
$p=1$& &Dirichlet & Robin/Neumann \\ \noalign{\smallskip}\hline\noalign{\smallskip}
$K_{n,i,h}$ & $\partial u(t_n)$ & - & - \\ \noalign{\smallskip}\hline\noalign{\smallskip}
 \multirow{2}{1cm}{$U_h^{n+1}$} &  $\partial u(t_n)$ & - & - \\
\cline{2-4} & $\partial A u(t_n)/ \partial \bar{F}_{n,i}$& - & $\displaystyle{O(\|e_{n,h}\|_h)}$  \\ \noalign{\smallskip}\hline
\end{tabular}
\end{table}

\begin{table}
\caption{Errors which are committed at each step when approximating the corresponding  simplifications of the boundary terms in Tables \ref{t1} and \ref{t2} to achieve local order $3$}
\label{t4}
\begin{tabular}{llcc}
\hline\noalign{\smallskip}
$p=2$& &Dirichlet & Robin/Neumann \\ \noalign{\smallskip}\hline\noalign{\smallskip}
\multirow{2}{1cm}{$K_{n,i,h}$} & $\partial u(t_n)$ & - & - \\
\cline{2-4} & $\partial A u(t_n)/ \partial \bar{F}_{n,i}$& - & $O(\|e_{n,h}\|_h)$ \\ \noalign{\smallskip}\hline\noalign{\smallskip}
\multirow{4}{1cm}{$U_h^{n+1}$} &  $\partial u(t_n)$ & - & - \\
\cline{2-4} & $\partial A u(t_n)$& - & $O(\|e_{n,h}\|_h)$ \\
\cline{2-4} &  $\partial \bar{F}_{n,i}$ & - & $O(\|e_{n,h}\|_h+k \mu_{k,1})$ \\
\cline{2-4} & $\partial A^2 u(t_n)/ \partial A \bar{F}_{n,i}$& $O(\nu_h+\frac{\|e_{n,h}\|_h}{h^\gamma})$ & $\displaystyle{O(\mu_{k,1}+\frac{\|e_{n,h}\|_h}{k}+\nu_h+\frac{\|e_{n,h}\|_h}{h^\gamma})}$ \\ \noalign{\smallskip}\hline
\end{tabular}
\end{table}

\begin{table}
\caption{Errors which are committed at each step when approximating the corresponding simplifications of the boundary terms in Tables \ref{t1} and \ref{t2}
  to achieve local order $4$}
\label{t5}
\begin{tabular}{llcc}
\hline\noalign{\smallskip}
$p=3$& &Dirichlet & Robin/Neumann \\ \noalign{\smallskip}\hline\noalign{\smallskip}
\multirow{4}{1cm}{$K_{n,i,h}$} & $\partial u(t_n)$ & - & - \\
\cline{2-4} & $\partial A u(t_n)$& - & $O(\|e_{n,h}\|_h)$ \\
\cline{2-4} & $\partial \bar{F}_{n,i}$ & - & $O(\|e_{n,h}\|_h+k \mu_{k,1})$ \\
\cline{2-4} & $\partial A^2 u(t_n)/ \partial A \bar{F}_{n,i}$ & $O(\nu_h+\frac{\|e_{n,h}\|_h}{h^\gamma})$ & $O(\mu_{k,1}+\frac{\|e_{n,h}\|_h}{k}+\nu_h+\frac{\|e_{n,h}\|_h}{h^\gamma})$ \\ \noalign{\smallskip}\hline\noalign{\smallskip}
\multirow{6}{1cm}{$U_h^{n+1}$} &  $\partial u(t_n)$ & - & - \\
\cline{2-4} & $\partial A u(t_n)$& - & $O(\|e_{n,h}\|_h)$ \\
\cline{2-4} &  $\partial \bar{F}_{n,i}$ & $O(k^2 \nu_h +\frac{k^2}{h^\gamma} \|e_{n,h}\|_h)$ & $O(\|e_{n,h}\|_h+k \mu_{k,1}+k^2 \mu_{k,2}+k^2 \nu_h+ \frac{k^2}{h^\gamma}\|e_{n,h}\|_h$)  \\
\cline{2-4} & $\partial A^2 u(t_n)$& $O(\nu_h+\frac{\|e_{n,h}\|_h}{h^\gamma})$ & $O(\mu_{k,1}+\frac{\|e_{n,h}\|_h}{k}+\nu_h+\frac{\|e_{n,h}\|_h}{h^\gamma})$ \\
\cline{2-4} &  $\partial A \bar{F}_{n,i}$ & $O(\nu_h+\frac{\|e_{n,h}\|_h}{h^\gamma}+\frac{k \mu_{k,1}}{h^\gamma})$ &  $O( \mu_{k,1}+k \mu_{k,2}+\frac{\|e_{n,h}\|_h}{k}+\nu_h+\frac{\|e_{n,h}\|_h}{h^\gamma})$\\
\cline{2-4} & $\partial A^3 u(t_n) / \partial A^2 \bar{F}_{n,i}$ & $O(\nu_h+\frac{\|e_{n,h}\|_h}{k h^\gamma}+\frac{ \mu_{k,1}}{h^\gamma})$ & $O( \mu_{k,1}+ \mu_{k,2}+\frac{\|e_{n,h}\|_h}{k^2}+\nu_h+\frac{e_{n,h}}{k h^\gamma})$ \\ \noalign{\smallskip}\hline
\end{tabular}

\end{table}

\subsection{Global error of the full discretization}

In this subsection we show how the global error $e_{n,h}$ behaves with our suggestion for the implementation of EERK methods, which is based in formulas (\ref{Knih})-(\ref{Uhn}), after approximating the boundaries of those formulas as suggested in the previous subsections. For that, the classical argument which relates local error to global error can be performed, but the error coming from the approximation of the simplified boundaries must also be considered. Besides, under certain assumptions related to more regularity and parabolicity, a summation-by-parts argument can be applied and therefore, the order can increase.

\begin{theorem}
Let us assume the first set of hypotheses of Theorem \ref{thlocs} and also, just for $p=2,3,$ that, for a certain constant $C$,
\begin{eqnarray}
\frac{k}{h^\gamma} \le C,
\label{cfl}
\end{eqnarray}
where $\gamma$ is the order of the space derivative which must be approximated through numerical differentiation to calculate the necessary boundaries of the suggested method. Then,
it happens that, under Dirichlet boundary conditions,
\begin{list}{$\bullet$}{}
\item For $p=1$,  $\|e_{n,h}\|_h=O(k+\varepsilon_h)$,
\item For $p=2$, $\|e_{n,h}\|_h=O(k^2+k \nu_h+\varepsilon_h)$,
\item For $p=3$, $\|e_{n,h}\|_h=O(k^3+k \nu_h+k \mu_{k,1}+\varepsilon_h)$,
\end{list}
and,  under R/N boundary conditions,
\begin{list}{$\bullet$}{}
\item For $p=1$,  $\|e_{n,h}\|_h=O(k+\varepsilon_h)$,
\item For $p=2$, $\|e_{n,h}\|_h=O(k^2+k \mu_{k,1}+k \nu_h+\varepsilon_h)$,
\item For $p=3$, $\|e_{n,h}\|_h=O(k^3+k \mu_{k,1}+k^2 \mu_{k,2}+k \nu_h+\varepsilon_h)$.
\end{list}
Assuming also the second set of hypotheses of the same theorem, that the following condition holds for a constant $C$ which is independent of $k$ and $h$
\begin{eqnarray}
\| k A_{h,0} \sum_{r=1}^{n-1} e^{r k A_{h,0}}\|_h \le C, \quad 0 \le nk \le T,
\label{parabol}
\end{eqnarray}
and that, for $t\in [0,T]$, $\dot{u}(t)\in D(A^{p+1})$ and
\begin{eqnarray}
A^l \dot{u}(t) \in Z, \quad l=0,1,\dots,p+1,
\label{regflesp1}
\end{eqnarray}
it happens that, under Dirichlet boundary conditions,
\begin{list}{$\bullet$}{}
\item For $p=1$,  $\|e_{n,h}\|_h=O(k^2+k \varepsilon_h+\eta_h)$,
\item For $p=2$, $\|e_{n,h}\|_h=O(k^3+k \nu_h+k \varepsilon_h+\eta_h)$,
\item For $p=3$, $\|e_{n,h}\|_h=O(k^4+k \nu_h+k \mu_{k,1}+k\varepsilon_h+\eta_h)$,
\end{list}
and, under R/N boundary conditions,
\begin{list}{$\bullet$}{}
\item For $p=1$,  $\|e_{n,h}\|_h=O(k^2+k \varepsilon_h+\eta_h)$,
\item For $p=2$, $\|e_{n,h}\|_h=O(k^3+k \mu_{k,1}+k \nu_h+k \varepsilon_h+\eta_h)$,
\item For $p=3$, $\|e_{n,h}\|_h=O(k^4+k \mu_{k,1}+k^2 \mu_{k,2}+k \nu_h+k\varepsilon_h+\eta_h)$.
\end{list}
\label{thglob}
\end{theorem}

Notice that condition (\ref{parabol}) was also used in \cite{acrnl} and was proved in \cite{HO0} for analytic semigroups, covering the case in which the linear operator in (\ref{laibvp}) corresponds to a parabolic problem.

On the other hand, we remark that, for $p=2,3$, a CFL condition (\ref{cfl}) which, in any case, is much less restrictive than that needed with explicit Runge-Kutta methods is required.

\section{Numerical experiments}

In this section, we have made some numerical experiments in order to corroborate the  results of the previous sections.

\subsection{One-dimensional problem}

\begin{table}
\caption{Local and global error when integrating problem (\ref{p1})  with nonvanishing boundary conditions (\ref{dir}) with method (\ref{rk2}) without avoding order reduction.}
\label{t6}
\begin{tabular}{ccccc} \hline\noalign{\smallskip}
k & 1/20 & 1/40 & 1/80 & 1/160 \\ \noalign{\smallskip}\hline\noalign{\smallskip}
Local error & 2.1370e-2 & 1.0531e-2 & 5.2115e-3 & 2.5816e-3 \\
Order & &1.02 & 1.01& 1.01 \\
Global error & 2.2693e-2 & 1.1280e-2 & 5.6078e-3 & 2.7848e-3 \\
Order & & 1.01 & 1.01 & 1.01\\ \noalign{\smallskip}\hline
\end{tabular}
\end{table}

\begin{table}
\caption{Local and global error when integrating  problem (\ref{p1})  with  nonvanishing boundary conditions (\ref{dn}) with method (\ref{rk2}) without avoding order reduction.}
\label{t6dn}
\begin{tabular}{ccccc} \hline\noalign{\smallskip}
k & 1/20 & 1/40 & 1/80 & 1/160 \\ \noalign{\smallskip}\hline\noalign{\smallskip}
Local error & 1.1214e-3 & 3.9698e-4 & 1.4089e-4 & 4.9963-5 \\
Order & & 1.50 & 1.49 & 1.50 \\
Global error & 2.0719e-2 & 1.0369e-2 & 5.1725e-3 & 2.5729e-3 \\
Order & & 0.99  & 1.00 & 1.01 \\ \noalign{\smallskip}\hline
\end{tabular}

\end{table}

\begin{table}
\caption{Local and global error when integrating problem (\ref{p1}) with  nonvanishing boundary conditions (\ref{dir}) with method (\ref{rk2}) using $p=1$ in (\ref{Knih})-(\ref{Uhn}) and the simplification for $\partial \bar{F}_{n,i}$ in $U_h^{n+1}$ in Table \ref{t2} for that value of $p$.}
\label{t7}
\begin{tabular}{ccccc} \hline\noalign{\smallskip}
k & 1/20 & 1/40 & 1/80 & 1/160 \\ \noalign{\smallskip}\hline\noalign{\smallskip}
Local error & 1.2438e-3 & 3.1031e-4 & 7.7348e-5 & 1.9256e-5 \\
Order & & 2.00 & 2.00& 2.01\\
Global error & 7.5946e-4 & 1.9141e-4 & 4.8515e-5 & 1.2292e-5 \\
Order & & 1.99 & 1.98 & 1.98\\ \noalign{\smallskip}\hline
\end{tabular}
\end{table}

\begin{table}
\caption{Local and global error when integrating problem (\ref{p1}) with nonvanishing boundary conditions (\ref{dn}) with method (\ref{rk2}) using $p=1$ in (\ref{Knih})-(\ref{Uhn}) and the simplification for $\partial \bar{F}_{n,i}$ in $U_h^{n+1}$ in Table \ref{t2} for that value of $p$.}
\label{t7dn}
\begin{tabular}{ccccc} \hline\noalign{\smallskip}
k & 1/20 & 1/40 & 1/80 & 1/160 \\ \noalign{\smallskip}\hline\noalign{\smallskip}
Local error &  1.2438e-3 & 3.1031e-4 & 7.7348e-5 & 1.9256e-5\\
Order & & 2.00 & 2.00 & 2.01\\
Global error & 8.1236e-4 & 2.0456e-4 & 5.1368e-5 & 1.2857e-5 \\
Order & & 1.99 & 1.99 & 2.00 \\ \noalign{\smallskip}\hline
\end{tabular}
\end{table}

\begin{table}
\caption{Local and global error when integrating problem (\ref{p1})  with nonvanishing boundary conditions (\ref{dir}) with method (\ref{rk2}) using $p=2$ in (\ref{Knih})-(\ref{Uhn}) and the simplification for the boundaries in Tables \ref{t1} and \ref{t2} for that value of $p$. }
\label{t8}
\begin{tabular}{ccccc} \hline\noalign{\smallskip}
k & 1/20 & 1/40 & 1/80 & 1/160 \\ \noalign{\smallskip}\hline\noalign{\smallskip}
Local error & 4.9174e-5 & 6.5311e-6 & 8.3433e-7 & 1.0548e-7 \\
Order & &2.91 & 2.97& 2.98\\
Global error & 4.2309e-5 & 1.0515e-5 & 2.6288e-6 & 6.5912e-7 \\
Order & & 2.01& 2.00& 2.00\\ \noalign{\smallskip}\hline
\end{tabular}
\end{table}

\begin{table}
\caption{Local and global error when integrating problem (\ref{p1})  with  nonvanishing boundary conditions (\ref{dn}) with method (\ref{rk2}) using $p=2$ in (\ref{Knih})-(\ref{Uhn}) and the simplification for the boundaries in Tables \ref{t1} and \ref{t2} for that value of $p$. }
\label{t8dn}
\begin{tabular}{ccccc} \hline\noalign{\smallskip}
k & 1/20 & 1/40 & 1/80 & 1/160 \\ \noalign{\smallskip}\hline\noalign{\smallskip}
Local error &  4.9884e-5 & 6.5394e-6 & 8.3434e-7 & 1.0548e-7\\
Order & & 2.93 & 2.97 & 2.98 \\
Global error & 2.0216e-4 & 4.9199e-5 & 1.2062e-5 & 2.8962e-6 \\
Order & & 2.04 & 2.03 & 2.06 \\ \noalign{\smallskip}\hline
\end{tabular}
\end{table}

\begin{figure}
\epsfig{file=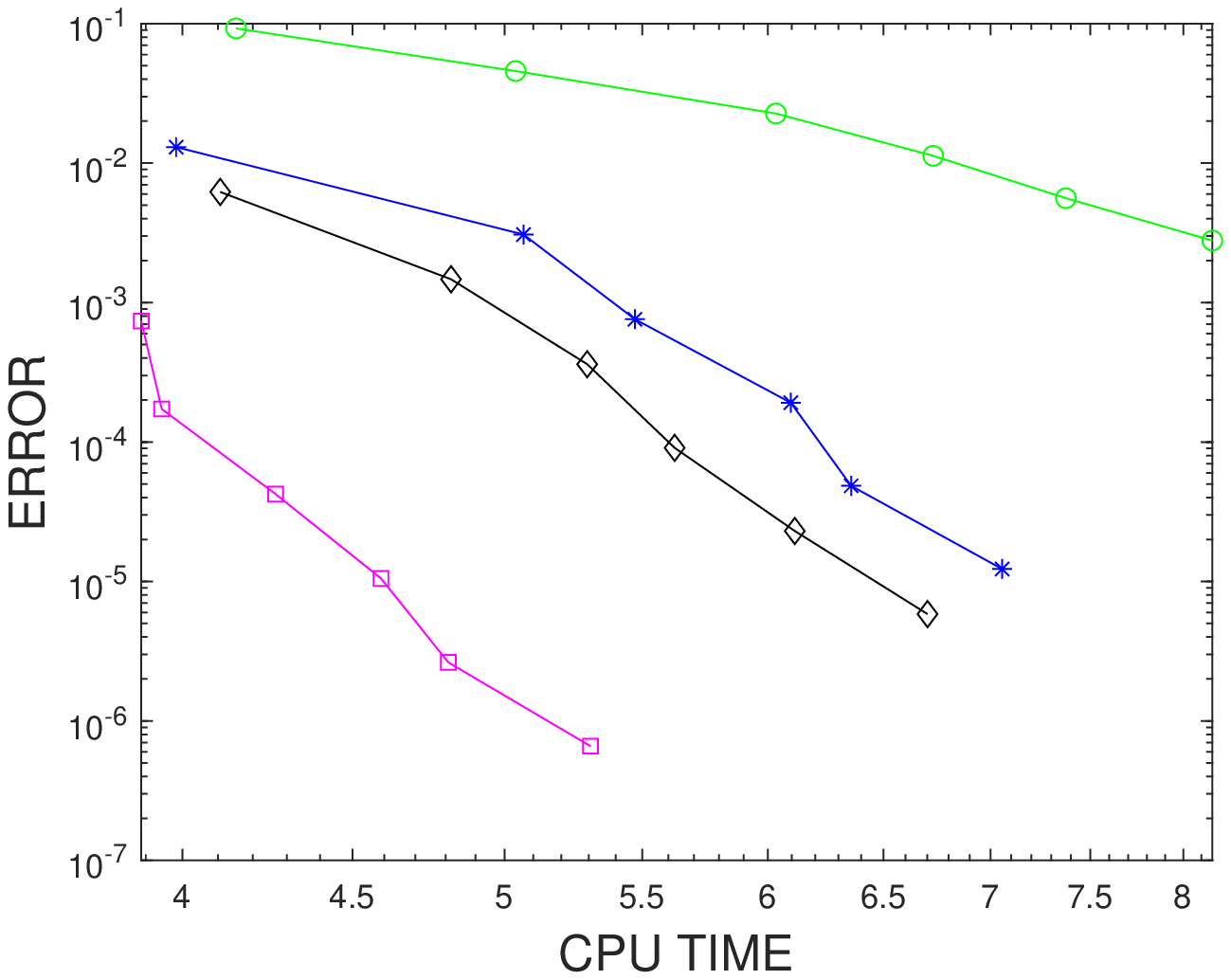,height=6in,width=5in}
\caption{Error against CPU time when integrating problem (\ref{p1}) with nonvanishing boundary conditions (\ref{dir}), using (\ref{rk2}) without avoiding order reduction
(green circles), the suggested technique (\ref{Knih1})-(\ref{Uhn}) with $p=1$ and the simplifications in Table \ref{t2} (blue asterisks), the suggested technique corresponding to $p=2$  and the simplifications in Tables \ref{t1} and \ref{t2} (\ref{rk2p2}) (magenta squares) and the method with stiff order $2$ (\ref{rk2b}) (black diamonds).}
\label{f1}
\end{figure}

\begin{figure}
\epsfig{file=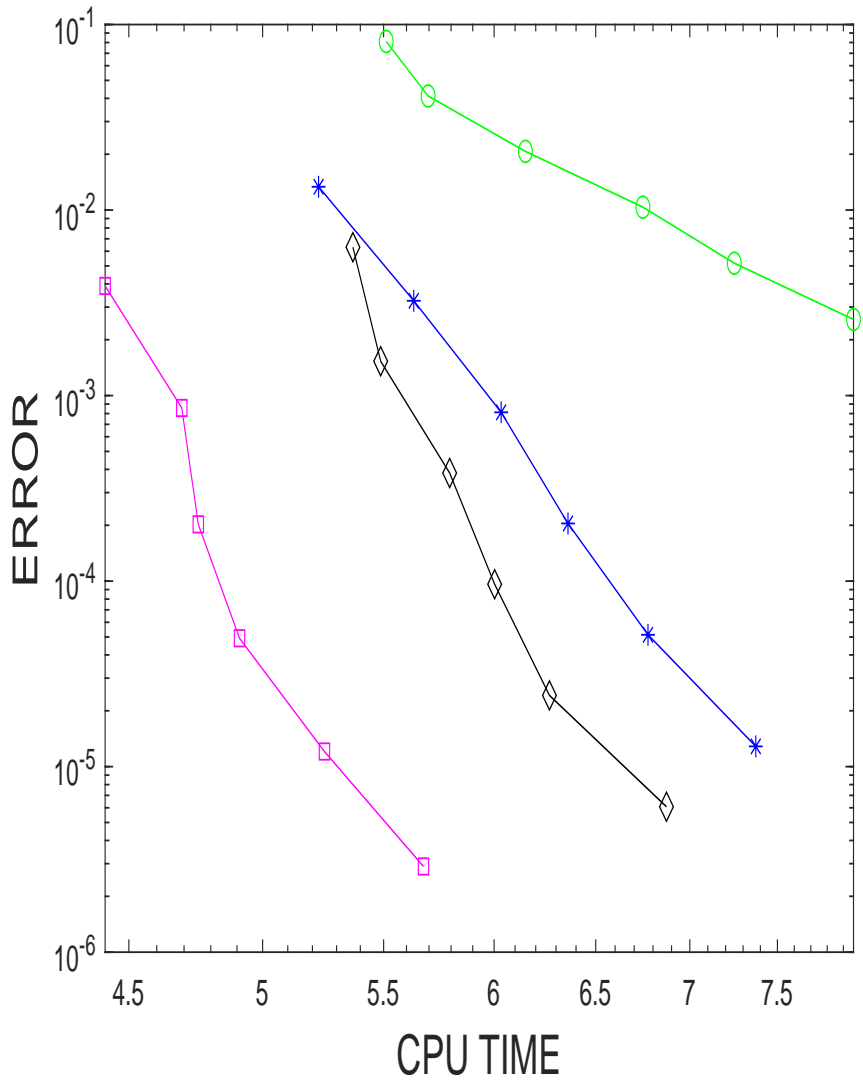,height=6in,width=5in}
\caption{Error against CPU time when integrating problem (\ref{p1}) with nonvanishing boundary conditions (\ref{dn}), using (\ref{rk2}) without avoiding order reduction
(green circles), the suggested technique (\ref{Knih1})-(\ref{Uhn}) with $p=1$ and the simplifications in Table \ref{t2} (blue asterisks), the suggested technique corresponding to $p=2$  and the simplifications in Tables \ref{t1} and \ref{t2} (\ref{rk2p2}) (magenta squares) and the method with stiff order $2$ (\ref{rk2b}) (black diamonds).}
\label{f1dn}
\end{figure}

In the first place, we have considered the one-dimensional problem
\begin{eqnarray}
u_t=u_{xx}+u^2 +h(x,t), \quad x \in [0,1], t\in [0,1],
\label{p1}
\end{eqnarray}
with function $h$ and initial and boundary conditions so that the exact solution is $u(x,t)=\cos(x+t)$. On the one hand, we have considered  Dirichlet boundary conditions at both sides
\begin{eqnarray}
u(0,t)=g_0(t), \quad u(1,t)=g_1(t),
\label{dir}
\end{eqnarray}
and, on the other hand, Dirichlet at $x=0$ and Neumann at $x=1$
\begin{eqnarray}
u(0,t)=g_0(t), \quad u_x(1,t)=g_1(t).
\label{dn}
\end{eqnarray}
In such a way, (A1)-(A8) are clearly satisfied with $X=C[0,1]$, $\|\cdot \|$ the maximum norm (see the example in page 231 in \cite{acrnl}), $\phi(u)=u^2$, $h(t)=-\sin(x+t)+\cos(x+t)-\cos^2(x+t)$, $A$ the second-order differential operator in space, $\bar{m}(A)=2$ in (A6) and $$\|v\|_l=\max(\|v\|_{\infty},\|v'\|_{\infty},\dots,\|v^{(2l)}\|_{\infty}), \quad \mbox{ for }v \in D(A^l)=C^{2l}([0,1]),$$
in (A7) and (A8).

Then, for the space discretization, we have taken the standard second-order difference scheme, for which $D_h\equiv 0$ in (\ref{spacediscr}) and,  in the case of two  Dirichlet boundary conditions,
$$A_{h,0}=\mbox{tridiag}(1,-2,1)/h^2, \quad C_h[g_0(t),g_1(t)]^T =[g_0(t),0,\dots,0,g_1(t)]^T/h^2,$$
while, in the Dirichlet and Neumann one, $A_{h,0}$ has one more row and column since the value at the final node is also numerically approximated. More precisely, by considering a centered difference in that last node to approximate the Neumann condition,
$$
A_{h,0}=\frac{1}{h^2} \left( \begin{array}{cccccc} -2 & 1 & 0 & \dots & \dots & 0 \\ 1 & -2 & 1 & 0 & \dots & \vdots \\ 0 & 1  & -2 & 1 & \ddots & \vdots \\
\vdots & 0 & \ddots & \ddots & \ddots & 0 \\ \vdots & & \ddots &1 & -2& 1 \\ 0 & \dots & \dots &0 &2 & -2 \end{array}\right), \quad C_h \left( \begin{array}{c} g_0(t) \\ g_1(t) \end{array} \right)=\left( \begin{array}{c} \frac{1}{h^2} g_0(t) \\ 0 \\ \vdots \\\frac{2}{h}g_1(t) \end{array} \right).
$$
It can be proved (see \cite{acrnl} and \cite{S}) that (H1)-(H3) are satisfied for these space discretizations when $Z=C^4([0,1])$, $\varepsilon_h=O(h^2)$ for the completely Dirichlet case and $\varepsilon_h=O(h)$ for the Dirichlet and Neumann one, $\eta_h=O(h^2)$ and $\bar{m}(Z)=4$.

As for the time integration, we have considered  the EERK method
\begin{eqnarray}
\begin{array}{c|cc} 0 &  & \\ \frac{1}{2} & \frac{1}{2}\varphi_{1,2} &  \\  \hline & 0 & \varphi_1 \end{array}, \label{rk2}
\end{eqnarray}
which was proved to lead to order $\ge 1$ when applied to problems with vanishing boundary conditions \cite{HO}. (We notice that $\varphi_{1,2}(z)=\varphi_1(z/2)$ and that, with the notation in (\ref{eerk}), $a_{21}(z)=\varphi_1(z/2)/2$, $b_1(z)=0$ and $b_2(z)=\varphi_1(z)$). When this method is applied directly to
\begin{eqnarray}
\dot{U}_h(t)&=& A_{h,0} U_h(t)+C_h g(t)+U_h.^2 +P_h h(t), \nonumber \\
U_h(0)&=& P_h u_0, \label{waor}
\end{eqnarray}
where $U_h.^2$ denotes the vector of the square of components of $U_h$. We are in fact taking $M=A_{h,0}$ and $F(t,U)=C_h g(t)+U_h.^2+P_h h(t)$ in (\ref{etapas})-(\ref{eerk}), and thus local and global time order around 1 can be observed with Dirichlet boundary conditions (\ref{dir}), as Table 6 shows when considering $h=1/1000$, so that the error in space is negligible. As for Dirichlet and Neumann boundary conditions case (\ref{dn}), a bit higher local order can be observed in Table \ref{t6dn}, although the global order continues to be around 1.

As this method has nonstiff order 2, the results in this paper can be applied for both $p=1$ and $p=2$. More precisely, when implementing the corresponding formulas in (\ref{Knih}) and (\ref{Uhn}) with $p=1$ and the simplification for $\partial \bar{F}_{n,i}$ in $U_h^{n+1}$ in Table \ref{t2} for that value of $p$, Theorem \ref{thlocs} assures local order 2 in the timestepsize when the error in space is small enough and Theorem \ref{thglob} global order 2. When implementing the same formulas with $p=2$ using Tables \ref{t1} and \ref{t2} for the simplification of some boundaries, Theorem \ref{thlocs} assures local order 3 and Theorem \ref{thglob} global order 2, in principle under a CFL condition of the form (\ref{cfl}) with $\gamma=1$ due to Dirichlet boundary conditions. We remind that, as stated in Section 4.5, this CFL condition comes from the numerical approximation of $\partial A^2 u(t_n)$ and $\partial A f(t_n,u(t_n))$ through numerical differentiation in space. However, it happens that, for this particular method, this numerical differentiation is not required. This is due to the fact that, for $p=2$, (\ref{Knih})-(\ref{Uhn}) read like this after  the simplification on the boundaries in Tables \ref{t1} and \ref{t2}:
\begin{eqnarray}
K_{n,1,h}&=&U_h^n, \nonumber \\
K_{n,2,h}&=& e^{\frac{k}{2}A_{h,0}}U_h^n+\frac{k}{2} \varphi_1(\frac{k}{2} A_{h,0})C_h \partial u(t_n)+\frac{k^2}{4} \varphi_2(\frac{k}{2} A_{h,0})C_h \partial A u(t_n), \nonumber \\
&&+\frac{k}{2}[\varphi_1(\frac{k}{2} A_{h,0})F_{n,1,h}+\frac{k}{2}\varphi_2(\frac{k}{2} A_{h,0})C_h \partial f(t_n,u(t_n))], \nonumber \\
U_h^{n+1}&=&e^{k A_{h,0}}U_h^n+k \varphi_1(k A_{h,0}) C_h \partial u(t_n)+k^2 \varphi_2(k A_{h,0}) C_h \partial A u(t_n)+ k^3 \varphi_3(k A_{h,0}) C_h \partial A^2 u(t_n)
\nonumber \\
&&+k[\varphi_1(k A_{h,0})F_{n,2,h}+k\varphi_2(k A_{h,0})C_h \partial f(t_n+\frac{k}{2},u(t_n)+\frac{k}{2}\dot{u}(t_n)) \nonumber \\
&&\hspace{1cm}+k^2 \varphi_3(k A_{h,0})C_h \partial A f(t_n, u(t_n))].
\label{rk2p2}
\end{eqnarray}
Then, because of (\ref{laibvp}),
$$A^2 u(t_n)+ A f(t_n,u(t_n))=A \dot{u}(t_n)=\ddot{u}(t_n)-f_t(t_n,u(t_n))-f_u(t_n,u(t_n))\dot{u}(t_n),$$
which means that the whole term which multiplies $k^3 \varphi_3(k A_{h,0})$ in (\ref{rk2p2}) can be exactly calculated in terms of the data of the problem, when considering Dirichlet boundary conditions on both sides. Therefore, for this particular problem and time integrator, numerical differentiation is not required and condition (\ref{cfl}) is not necessary. On the other hand, when considering the Neumann boundary condition on one side, no space numerical differentiation is either required, but time numerical integration is needed indeed. For that, we have used a second-order Taylor expansion for the first time step and then, a 2-BDF formula for the next ones.

Tables \ref{t7},  \ref{t7dn}, \ref{t8} and 11 show that the expected numerical order is in fact achieved. When implementing the different techniques which have been described here, putting together all the terms which multiply each $\varphi_j$ and using Krylov subroutines \cite{niesen}, we can see that procedure (\ref{rk2p2}) is the most efficient for both types of boundary conditions. (See Figures \ref{f1} and \ref{f1dn}, where it is clear that, for a same amount of CPU time, the technique with the smallest error is the one corresponding to the suggested technique with $p=2$. A thorough explanation for that is given in \cite{CR1}). Although it is not an aim of this paper, we have also compared in the same figures with another EERK method in \cite{HO},
\begin{eqnarray}
\begin{array}{c|cc} 0 &  & \\ \frac{1}{2} & \frac{1}{2}\varphi_{1,2} &  \\  \hline & \varphi_1-2 \varphi_2 & 2 \varphi_2 \end{array}, \label{rk2b}
\end{eqnarray}
which happens to have stiff order $2$ for problem (\ref{p1}). (Notice that this method and (\ref{rk2}) have the same underlying RK method.) When implementing this method by integrating directly (\ref{waor}) with the formulas in (\ref{eerk}), order $2$ is observed, but is less efficient than the procedure in (\ref{rk2p2}).

\subsection{Bidimensional problem}

\begin{table}
\caption{Local and global error when integrating Dirichlet problem (\ref{p2})  with nonvanishing boundary conditions with method (\ref{krogstad}) without avoding order reduction.}
\label{t9}
\begin{tabular}{ccccccc} \hline\noalign{\smallskip}
k & 1/8 & 1/16 & 1/32 & 1/64 & 1/128& 1/256 \\ \noalign{\smallskip}\hline\noalign{\smallskip}
Local error & 3.6983e-5 & 4.9435e-6 & 6.6113e-7 & 8.6583e-08 & 1.1152e-8 & 1.4159e-9 \\
Order &  & 2.90 & 2.90 & 2.93 & 2.96 & 2.98  \\
Global error & 3.4712e-5 & 4.4791e-6 & 5.9470e-7 & 7.7735e-8 & 1.0010e-8 & 1.2770e-9   \\
Order & & 2.95 & 2.91 & 2.94 & 2.96 & 2.97   \\ \noalign{\smallskip}\hline
\end{tabular}

\end{table}

\begin{table}[t]
\caption{Local and global error when integrating Dirichlet problem (\ref{p2})  with nonvanishing boundary conditions with method (\ref{krogstad}) using $p=3$ in (\ref{Knih})-(\ref{Uhn}) and the simplification for the boundaries in Tables \ref{t1} and \ref{t2} for that value of $p$.}
\label{t10}
\begin{tabular}{ccccccc} \hline\noalign{\smallskip}
k & 1/8 & 1/16 & 1/32 & 1/64 & 1/128& 1/256 \\ \noalign{\smallskip}\hline\noalign{\smallskip}
Local error & 8.2450e-5 & 4.9697e-6 & 3.0088e-7 & 1.8147e-8 & 1.0834e-9 & 6.3605e-11 \\
Order &  & 4.05 & 4.05 & 4.05 & 4.07 & 4.09 \\
Global error & 1.0842e-4 & 6.9957e-6 & 4.3923e-7 & 2.7503e-8 & 1.7108e-9 & 1.0553e-10 \\
Order & & 3.95 & 3.99 & 4.00 & 4.01 & 4.02 \\ \noalign{\smallskip}\hline
\end{tabular}
\end{table}

\begin{figure}
\epsfig{file=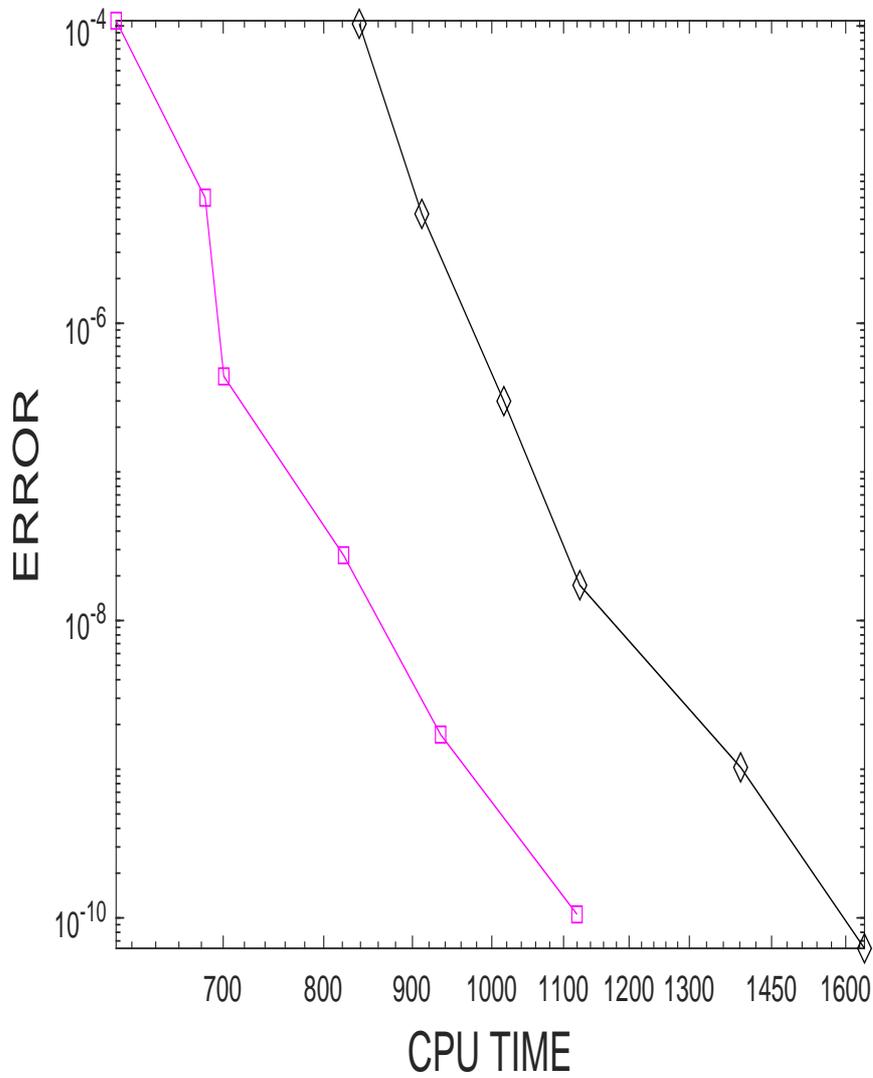,height=6in,width=5in}
\caption{Error against CPU time when integrating problem (\ref{p2}), using (\ref{krogstad}) with the suggested technique (\ref{Knih})-(\ref{Uhn}) with $p=3$ and the simplifications in Tables \ref{t1} and \ref{t2} (magenta squares), and the method in \cite{LO} with stiff order $4$ when applied directly to (\ref{sdd}) (black diamonds).}
\label{f2}
\end{figure}

In a second place, we have considered the bidimensional problem
\begin{eqnarray}
u_t=u_{xx}+u_{yy}+u^2+h(x,y,t), \quad x,y\in [0,1], \quad t\in [0,1], \label{p2}
\end{eqnarray}
with function $h$ and initial and Dirichlet boundary conditions such that the exact solution is $u(x,y,t)=\cos(t+x+y)$ and again (A1)-(A8) are satisfied. We notice that
$$h(x,y,t)=-\sin(t+x+y)+2 \cos(t+x+y)-\cos^2(t+x+y).$$
As space discretization, we have considered the $9$-point formula \cite{S}, which determines certain matrices $A_{h,0}$, $C_h$ and $D_h$ in (\ref{spacediscr}). It can be observed that (H1)-(H3) are again satisfied (see \cite{acr2,S}) for $Z=C^6 ([0,1]\times [0,1])$ with $\varepsilon_h=O(h^4)$, $\eta_h=O(h^4)$ and $\bar{m}(Z)=6$. As for the time integration, we consider Krogstad's method \cite{K}
\begin{eqnarray}
\begin{array}{c|cccc} 0 &  &  & & \\ \frac{1}{2} & \frac{1}{2}\varphi_{1,2} & & & \\ \frac{1}{2} & \frac{1}{2}\varphi_{1,3}-\varphi_{2,3}& \varphi_{2,3} &  & \\ 1 & \varphi_{1,4}-2 \varphi_{2,4}& 0 & 2\varphi_{2,4} & \\
 \hline & \varphi_1-3 \varphi_2+4 \varphi_3 & 2\varphi_2-4 \varphi_3& 2\varphi_2-4 \varphi_3 & 4\varphi_3-\varphi_2 \end{array}. \label{krogstad}
\end{eqnarray}
 This method was proved to lead to order $3$ when applied to stiff problems with vanishing boundary conditions in \cite{HO}, although its nonstiff order is well known to be 4 \cite{K}. We have implemented the method when integrating directly
\begin{eqnarray}
 \dot{U}_h(t)&=&A_{h,0} U_h(t)+C_h g(t)+U_h.^2+P_h h(t)+D_h[(g(t)^2+\partial h(t)-\dot{g}(t)], \nonumber \\
 U_h(0)&=&P_h u(0), \label{sdd}
\end{eqnarray}
through formulas (\ref{etapas})-(\ref{eerk}) and Krylov methods \cite{niesen} to calculate the exponential operators (We notice that, as the matrix $A_{h,0}$ is not so sparse in this case because it is just the product of the inverse of a sparse mass matrix times another sparse matrix, a modification of the subroutines in \cite{niesen} has been used following the lines of Section 7.1 in \cite{GG}, so as to take profit of the highly sparse structure of the underlying matrices.) Local and global order $3$ then turn up, as Table 9 shows for the space grid $h=1/160$. We have then applied the suggested formulas (\ref{Knih})-(\ref{Uhn}) with $p=3$ with this method, considering the simplifications in Tables 1 and 2. We notice that, for this particular method, no numerical differentiation is required for the stages since the terms in $\partial A^2 u(t_n)$ and $\partial A \bar{F}_{n,i}$ sum up to a term in $\partial A \dot{u}$, which can be calculated in terms of data because $A \dot{u}=\ddot{u}-\dot{h}-\phi'(u)\dot{u}$. However, for the approximation of $\partial A^2 u(t_n)/\partial \bar{F}_{n,i}$ in $U_h^{n+1}$, $u_x$ and $u_y$ are required at the boundary, for which we have used a 4th-order BDF formula in the sides of the square in which we could not calculate it directly in terms of $g$. As for the approximation of $\partial A^3 u(t_n)/\partial A \bar{F}_{n,i}/ \partial A^2 \bar{F}_{n,i}$, we required not only $u_x$, $u_y$, but also $\dot{u}_x$, $\dot{u}_y$ at the boundary. To calculate them, we have firstlyapproximated $\dot{u}$ in some points near the boundary, but at the interior of the domain. We have done that by using Taylor expansions of order $3$ near the initial condition, i.e.,
  $$
  u_t(x,y,lk) \approx u_t(x,y,0)+l k u_{tt}(x,y,0)+\frac{(lk)^2}{2} u_{ttt}(x,y,0), \quad l=0,1,2,
  $$
  where the terms above can be calculated in terms of $u_0$ by using equation (\ref{p2}).
  Then, after the third step, we have used the 3rd-order BDF formula
  $$
  u_t(x,y,lk) \approx \frac{1}{k}[\frac{11}{6}u(x,y,lk)-3 u(x,y,(l-1)k)+\frac{3}{2} u(x,y,(l-2)k)-\frac{1}{3} u(x,y,(l-3)k)].
  $$
  Once those approximations of $\dot{u}$ were calculated, we have approximated $\dot{u}_x$, $\dot{u}_y$ with a 4th-order BDF formula in space. In such a way, $\mu_{k,1}=O(k^3)$ and $\nu_h=O(h^4)$ in Theorem \ref{thglob}, so that the error in space is easy to be negligible and the fourth order in time can be seen. In particular, condition (\ref{cfl}) does not seem very restrictive since local and global order 4 turn up, as it was justified by Theorems \ref{thlocs} and \ref{thglob}, and Table 10 corroborates.

Again, although it is not an aim of this paper, we have made a comparison in CPU time of this method with another one which shows order 4 without avoiding order reduction. We firstly tested with method (5.19) in \cite{HO} because it was suggested there as a $4$th-order method for stiff problems with vanishing boundary conditions, but it happened to show just order $3$ with non-vanishing boundary values. Because of that, we decided to consider the $5$th-order method for stiff problems and vanishing boundary conditions which is constructed in \cite{LO}. This method is $4$th-order accurate when applied directly to the space discretized system (\ref{sdd}). The comparison in CPU time is shown in Figure \ref{f2}, where it is shown that the latter is about $1/2$ more expensive than the suggested modification of Krogstad's method when avoiding order reduction.

\section*{Acknowledgements}
This research has been supported by Ministerio de Ciencia e
Innovaci\'on and Regional Development European Funds through project PGC2018-101443-B-I00 and by Junta de Castilla y Le\'on and Feder through project VA169P20.

\section{Appendix}
\subsection{Proof of Theorem \ref{thloc}}
Notice that, for the first set of hypotheses, through its definition and using (\ref{phi0r})-(\ref{phijr}),
\begin{eqnarray}
\bar{u}_{n+1}&=&\phi_{0,p,u(t_n),u(t_n)}(k)+k\sum_{i=1}^s \sum_{l=1}^s \mu_{i,l} \phi_{l,p-1,\bar{F}_{n,i},\bar{F}_{n,i}}(k) \nonumber \\
&=&\sum_{l=0}^p \frac{k^l}{l!} A^l u(t_n)+k^{p+1} \varphi_{p+1}(k A_0) A^{p+1} u(t_n) \nonumber \\
&&+k \sum_{i=1}^s \sum_{l=1}^s \mu_{i,l}[\sum_{j=0}^{p-1} \frac{k^j}{(j+l)!} A^j \bar{F}_{n,i}+k^p \varphi_{p+l}(k A_0) A^p \bar{F}_{n,i}],
\label{ub}
\end{eqnarray}
where $\bar{F}_{n,i}=f(t_n+c_i k, \bar{K}_{n,i})$, with
\begin{eqnarray}
\bar{K}_{n,i}&=&\phi_{0,p-1,u(t_n),u(t_n)}(c_i k)+k\sum_{j=1}^{i-1} \sum_{l,r=1}^s \lambda_{i,j,l,r}\phi_{l,p-2,\bar{F}_{n,j},\bar{F}_{n,j}}(c_r k) \nonumber \\
&=&\sum_{l=0}^{p-1} \frac{(c_i k)^l}{l!} A^l u(t_n)+(c_i k)^p \varphi_p(c_i k A_0) A^p u(t_n) \nonumber \\
&&+k \sum_{j=1}^{i-1} \sum_{l,r=1}^s \lambda_{i,j,l,r}[\sum_{ll=0}^{p-2} \frac{(c_r k)^{ll}}{(ll+l)!} A^{ll} \bar{F}_{n,j}+(c_r k)^{p-1} \varphi_{p-1+l}(c_r k A_0)A^{p-1} \bar{F}_{n,j}]. \label{bkni}
\end{eqnarray}
 In the following, we will inductively prove that $\bar{K}_{n,i}\in D(A^p)$ and $A^p \bar{K}_{n,i}$ is uniformly bounded in $k$ for small enough $k$, so that (\ref{ub}) and (\ref{bkni}) have perfect sense because of (A6) and the assumed regularity on $\phi$ and $h$. The result comes from Lemma 3.1 in \cite{acr2} because, according to it,
\begin{eqnarray}
A \phi_{0,p-1,u(t_n),u(t_n)}(\tau)=\sum_{l=0}^{p-2} \frac{\tau^l}{l!}A^{l+1} u(t_n)+\tau^{p-1} \varphi_{p-1}(\tau A_0) A^p u(t_n), \label{for1}
\end{eqnarray}
and, arguing in the same way,
\begin{eqnarray}
A^2 \phi_{0,p-1,u(t_n),u(t_n)}(\tau)&=&\sum_{l=0}^{p-3} \frac{\tau^l}{l!}A^{l+2} u(t_n)+\tau^{p-2} \varphi_{p-2}(\tau A_0) A^p u(t_n), \nonumber \\
\vdots &=& \vdots \nonumber \\
A^p \phi_{0,p-1,u(t_n),u(t_n)}(\tau)&=& \varphi_0(\tau A_0) A^p u(t_n). \label{for2}
\end{eqnarray}
From this, $\bar{K}_{n,1}\in D(A^p)$, which implies that $\bar{F}_{n,1}\in D(A^p)$. As for $\bar{K}_{n,2}$, in a similar way we notice that
\begin{eqnarray}
A \phi_{j,p-2,\bar{F}_{n,1},\bar{F}_{n,1}}(\tau)&=&\sum_{l=0}^{p-3} \frac{\tau^l}{(l+j)!}A^{l+1} \bar{F}_{n,1}+\tau^{p-2} \varphi_{j+p-2}(\tau A_0) A^{p-1} \bar{F}_{n,1}, \nonumber \\
\vdots &=& \vdots \nonumber \\
A^{p-1} \phi_{j,p-2,\bar{F}_{n,1},\bar{F}_{n,1}}(\tau)&=&\varphi_j(\tau A_0) A^{p-1} \bar{F}_{n,1}, \nonumber \\
A^p \phi_{j,p-2,\bar{F}_{n,1},\bar{F}_{n,1}}(\tau)&=& \frac{1}{\tau}[\varphi_{j+1}(\tau A_0)-\frac{1}{(j+1)!}I] A^{p-1}  \bar{F}_{n,1}, \quad \mbox{ for } \tau \neq 0,\label{for3}
\end{eqnarray}
and therefore,
$$k A^p \phi_{j,p-2,\bar{F}_{n,1},\bar{F}_{n,1}}(c_r k)=\bigg{\{} \begin{array}{ll}
\frac{k}{j!} A^p \bar{F}_{n,1}, & \mbox{ if }c_r =0, \\
\frac{1}{c_r}[\varphi_{j+1}(c_r k A_0)-\frac{1}{(j+1)!}I] A^{p-1}  \bar{F}_{n,1}, & \mbox{ if }c_r \neq 0, \end{array} $$
from what $\bar{K}_{n,2}\in D(A^p)$ and $A^p \bar{K}_{n,2}$ is uniformly bounded on $k$ for small enough $k$, which implies the same fact for $\bar{F}_{n,2}$.
Proceeding inductively with the rest of stages, it is justified that (\ref{ub}) and (\ref{bkni}) are true.

Then, developing in Taylor series (which can be done because of (A7) and the assumed regularity), the terms in $k^j$ ($0\le j \le p-1$)  in (\ref{bkni})  have the same form as those in  (\ref{Kbord}) but substituting $M$ by the operator $A$, $U(t_n)$ by $u(t_n)$ and $F_{*,n}$ by $f_*(t_n,u(t_n))$. Moreover, the terms in $k^j (0\le j \le p)$ in (\ref{ub}) have the same form as those in (\ref{Ubord}) with the same type of substitution. We notice that the residue terms in the Taylor expansions around $u(t_n)$ of $A^j \phi(\bar{K}_{n,i})$ ($j=0,\dots,p-1$) in (\ref{ub}) and (\ref{bkni}) correspond to
$$
\frac{1}{(p-j)!}\int_0^1 A^j \big[ \phi^{(p-j)}\big((1-\lambda)u(t_n)+\lambda \bar{K}_{n,i}\big)(\bar{K}_{n,i}-u(t_n))^{p-j}\big] d \lambda.
$$
Using that $(1-\lambda)u(t_n)+\lambda \bar{K}_{n,i} \in D(A^j)$ for $\lambda\in [0,1]$ and that $\|\bar{K}_{n,i}-u(t_n)\|_j=O(k)$, as justified through (\ref{for1}),(\ref{for2}) and (\ref{for3}), from (A7) it follows that these residues are $O(k^{p-j})$ whenever $\phi \in C^{p-j+j\bar{m}(A)}$. We notice that this regularity of $\phi$ for $j=0,\dots,p-1$ is satisfied under the assumed regularity. On the other hand, the factor $k^{j+1}$  which is multiplying $A^j \bar{F}_{n,i}$ in (\ref{ub}) and (\ref{bkni}) leads to the fact that all the residues are $O(k^{p+1})$.

  At the same time, the same happens for the Taylor series of $u(t_{n+1})$, which terms in $k^j$ ($0\le j \le p$) also coincide with those in (\ref{taylorUord}) except for the previous substitution of terms. As those terms coincided for an ordinary differential system, they also coincide now till the terms in which we are interested.

On the other hand, if in formulas (\ref{etsro})-(\ref{usro}), $p$ is substituted by $\hat{p}< p$, the same argument leads to the fact that $\rho_n=O(k^{\hat{p}+1})$. Moreover, when applying $A_0^{-1}$, the first term in $k^{\hat{p}+1}$ in (\ref{ub}) would be
\begin{eqnarray}
A_0^{-1} k^{\hat{p}+1} \varphi_{\hat{p}+1}(k A_0) A^{\hat{p}+1} u(t_n)&=&A_0^{-1} k^{\hat{p}+1} [ k A_0 \varphi_{\hat{p}+2}(k A_0)+\frac{1}{(\hat{p}+1)!}] A^{\hat{p}+1} u(t_n) \nonumber \\
&=& k^{\hat{p}+2} \varphi_{\hat{p}+2}(k A_0) A^{\hat{p}+1} u(t_n)+\frac{k^{\hat{p}+1}}{(\hat{p}+1)!}A_0^{-1} A^{\hat{p}+1} u(t_n), \nonumber
\end{eqnarray}
where the formula (\ref{recurf}) has been used in the first equality.
We notice that the first term is $O(k^{\hat{p}+2})$ and the second one coincides with one corresponding to the power $k^{\hat{p}+1}$ in $A_0^{-1} u(t_{n+1})$.
In a similar way, the other terms in  $k^{\hat{p}+1}$ in $A_0^{-1}$ times (\ref{ub}) would be
\begin{eqnarray}
\mu_{i,l} A_0^{-1}  k^{\hat{p}+1} \varphi_{\hat{p}+l}(k A_0) A^{\hat{p}} \bar{F}_{n,i}^{\hat{p}}=\mu_{i,l} \bigg[k^{\hat{p}+2} \varphi_{\hat{p}+l+1}(k A_0) A^{\hat{p}} \bar{F}_{n,i}^{\hat{p}}+\frac{k^{\hat{p}+1}}{(\hat{p}+l)!}A_0^{-1} A^{\hat{p}}\bar{F}_{n,i}^{\hat{p}}\bigg], \nonumber
\end{eqnarray}
where $\bar{F}_{n,i}^{\hat{p}}=f(t_n+c_i k, \bar{K}_{n,i}^{\hat{p}})$ with $\bar{K}_{n,i}^{\hat{p}}$ that in (\ref{bkni}) with $p$ substituted by $\hat{p}$. We notice again that the first term is $O(k^{\hat{p}+2})$ and the second coincides with one of the others corresponding to the power $k^{\hat{p}+1}$ in $A_0^{-1} u(t_{n+1})$ considering that the difference between $\bar{K}_{n,i}$ and $\bar{K}_{n,i}^{\hat{p}}$ is at least $O(k)$ (due to $\hat{p}\ge 1$), hypothesis (A4) and taking (\ref{Ubord}) into account and the fact that the terms of that type there correspond to the same terms in $U(t_{n+1})$.

\subsection{Proof of Theorem \ref{thlocfd}}

We notice that $\rho_{n,h}$ can be decomposed as
$$\rho_{n,h}=(\bar{U}_h^{n+1}-P_h \bar{u}_{n+1})+P_h \rho_n,$$
where $\rho_n$ is the local error associated to the time semidiscretization and can be bounded as Theorem \ref{thloc} states. Using (\ref{acotp}) and (H3d), it suffices to prove then that $\bar{U}_h^{n+1}-P_h \bar{u}_{n+1}=O(k \varepsilon_h)$. For that, we firstly take into account that
\begin{eqnarray}
\bar{U}_{h}^{n+1}-P_h \bar{u}_{n+1}&=&\phi_{0,p,P_h u(t_n),u(t_n)}^h (k)-P_h \phi_{0,p,u(t_n),u(t_n)}(k) \nonumber \\
&&+k \sum_{i=1}^s\sum_{l=1}^s \mu_{i,l}[\phi_{l,p-1,\bar{F}_{n,i,h}, \bar{F}_{n,i}}^h (k)-P_h \phi_{l,p-1,\bar{F}_{n,i},\bar{F}_{n,i}}(k)].
\label{f0}
\end{eqnarray}
By definition of $\phi_{0,r,P_h u(t_n),u(t_n)}^h$ and $\phi_{0,r,u(t_n),u(t_n)}$ and using (\ref{rh}), it follows that
\begin{eqnarray}
\frac{d}{d \tau}[\phi_{0,p,P_h u(t_n),u(t_n)}^h (\tau)-P_h \phi_{0,p,u(t_n),u(t_n)}(\tau)]&=& A_{h,0}[\phi_{0,p,P_h u(t_n),u(t_n)}^h (\tau)-P_h \phi_{0,p,u(t_n),u(t_n)}(\tau)]\nonumber \\
&&+A_{h,0}(P_h-R_h)\phi_{0,p,u(t_n),u(t_n)}(\tau) \nonumber \\
\phi_{0,p,P_h u(t_n),u(t_n)}^h (0)-P_h \phi_{0,p,u(t_n),u(t_n)}(0)&=&0, \nonumber
\end{eqnarray}
which implies that
\begin{eqnarray}
\lefteqn{\phi_{0,p,P_h u(t_n),u(t_n)}^h (k)-P_h \phi_{0,p,u(t_n),u(t_n)} (k)} \nonumber \\
&=&\int_0^k e^{(k-s) A_{h,0}} A_{h,0}(P_h-R_h) \phi_{0,p,u(t_n),u(t_n)} (s) ds =O(k \varepsilon_h),
\label{fo1}
\end{eqnarray}
by using (\ref{consistency}) and that $\phi_{0,p,u(t_n),u(t_n)}\in C([0,T],Z)$ because of (\ref{phi0r}), (\ref{regfle}) and (H2a).
In a similar way,
\begin{eqnarray}
&&\frac{d}{d\tau}[\phi_{l,p-1,\bar{F}_{n,i,h},\bar{F}_{n,i}}^h (\tau)- P_h \phi_{l,p-1,\bar{F}_{n,i},\bar{F}_{n,i}}(\tau)]
  =(A_{h,0}-\frac{l}{\tau}I_h)[\phi_{l,p-1,\bar{F}_{n,i,h},\bar{F}_{n,i}}^h (\tau)- P_h \phi_{l,p-1,\bar{F}_{n,i},\bar{F}_{n,i}}(\tau)]\nonumber \\
&&\hspace{3cm}+\frac{1}{(l-1)! \tau}[\bar{F}_{n,i,h}-P_h \bar{F}_{n,i}]+A_{h,0}(P_h-R_h)\phi_{l,p-1,\bar{F}_{n,i},\bar{F}_{n,i}}(\tau) \nonumber \\
&&\phi_{l,p-1,\bar{F}_{n,i,h},\bar{F}_{n,i}}^h (0)-P_h \phi_{l,p-1,\bar{F}_{n,i},\bar{F}_{n,i}}(0)=\frac{1}{l!}[\bar{F}_{n,i,h}-P_h \bar{F}_{n,i}], \nonumber
\end{eqnarray}
where $I_h$ is the identity matrix with the same size as $A_{h,0}$.
Now, we take into account that the solution of this system is $V_h(\tau)+W_h(\tau)$, with $V_h(\tau)$ the solution of
\begin{eqnarray}
V_h'(\tau)&=& (A_{h,0}-\frac{l}{\tau}I_h)V_h(\tau)+A_{h,0}(P_h-R_h)\phi_{l,p-1,\bar{F}_{n,i},\bar{F}_{n,i}}(\tau) \nonumber \\
V_h(0)&=&0, \nonumber
\end{eqnarray}
and $W_h$ the solution of
\begin{eqnarray}
W_h'(\tau)&=& (A_{h,0}-\frac{l}{\tau}I_h)W_h(\tau)+\frac{1}{(l-1)! \tau}[\bar{F}_{n,i,h}-P_h \bar{F}_{n,i}] \nonumber \\
W_h(0)&=&\frac{1}{l!}[\bar{F}_{n,i,h}-P_h \bar{F}_{n,i}]. \nonumber
\end{eqnarray}

From this,
\begin{eqnarray}
V_h(k)&=& \int_0^k e^{\int_s^k (A_{h,0}-\frac{l}{\tau}I) d\tau} A_{h,0}(P_h-R_h)\phi_{l,p-1,\bar{F}_{n,i},\bar{F}_{n,i}}(s) ds \nonumber \\
&=& \int_0^k e^{(k-s)A_{h,0}}(\frac{s}{k})^l A_{h,0}(P_h-R_h)\phi_{l,p-1,\bar{F}_{n,i},\bar{F}_{n,i}}(s) ds \nonumber \\
&=&O(k \varepsilon_h), \label{fo2}
\end{eqnarray}
where again (\ref{consistency}) has been used noting that, inductively, $\phi_{l,p-1,\bar{F}_{n,i},\bar{F}_{n,i}}\in C([0,T],Z)$ because of (\ref{phijr}), (\ref{regfle}), the regularity of $\phi$, (\ref{for2}), (\ref{for3}) and (H2a). On the other hand, we notice that, using Lemma 6 in \cite{CaM},
\begin{eqnarray}
W_h (k)=\varphi_l(k A_{h,0})[\bar{F}_{n,i,h}-P_h \bar{F}_{n,i}].
\label{fo3}
\end{eqnarray}
As $\varphi_l(k A_{h,0})$ is bounded because of (H1b) and, because of (\ref{fandph}), $P_h \bar{F}_{n,i}=f(t_n+c_i k, P_h \bar{K}_{n,i})$, using (H3), what we need is to bound
\begin{eqnarray}
\bar{K}_{n,i,h}-P_h \bar{K}_{n,i}&=& \phi_{0,p-1,P_h u(t_n),u(t_n)}^h (c_i k)-P_h \phi_{0,p-1,u(t_n),u(t_n)}(c_i k) \nonumber \\
&&+k \sum_{j=1}^{i-1}\sum_{l,r=1}^s \lambda_{i,j,l,r}[\phi_{l,p-2,\bar{F}_{n,j,h}, \bar{F}_{n,j}}^h (c_r k)-P_h \phi_{l,p-2,\bar{F}_{n,j},\bar{F}_{n,j}}(c_r k)].
\nonumber
\end{eqnarray}
Using the same arguments as above, it can be proved, recursively in $i=1,\dots,s$, that $\bar{K}_{n,i,h}-P_h \bar{K}_{n,i}=O(k \varepsilon_h)$, from what, considering that in (\ref{fo3}) and inserting it in (\ref{f0}) together with (\ref{fo2}) and (\ref{fo1}), the first result follows.

As for the result concerning the second set of hypotheses, it suffices to decompose again
$$
A_{h,0}^{-1} \rho_{n,h}=A_{h,0}^{-1}(\bar{U}_{n+1,h}-P_h \bar{u}_{n+1})+A_{h,0}^{-1} P_h \rho_{n}.
$$
Then, in a similar way as above, but using now the second formula of (\ref{consistency}) when possible, it can be proved that $A_{h,0}^{-1}(\bar{U}_{n+1,h}-P_h \bar{u}_{n+1})=O(k \eta_h+k^2 \varepsilon_h)$. On the other hand, as $\partial A_0^{-1} \rho_{n+1}=0$, by using (\ref{rh}) and the fact that, at the next step, the boundary of the numerical solution will be considered as the exact, $A_{h,0} R_h (A_0^{-1} \rho_{n})=P_h \rho_{n}$, from what
$$
A_{h,0}^{-1} P_h \rho_n=P_h (A_0^{-1} \rho_{n})+(R_h-P_h)A_0^{-1} \rho_{n}.$$
Now, the first term on the right-hand side is $O(k^{\hat{p}+2})$ by using the second part of Theorem \ref{thloc} and (\ref{acotp}). As for the second term, it is $O(k \eta_h)$  because, for some $t^*\in(t_n,t_{n+1})$, $(u(t_{n+1})-u(t_n))/k= \dot{u}(t^*)=Au(t^*)+f(t^*, u(t^*))\in Z$ and $(\bar{u}_{n+1}-u(t_n))/k\in Z$, due to (\ref{regfle}) and (H2a).
Therefore,  $A_0^{-1} \rho_{n}/k \in Z$ and (\ref{consistency}) can be applied.

\subsection{Proof of Theorem \ref{thlocs}}

We firstly consider $p=1$ and notice that the same decomposition as in the proof of Theorem \ref{thlocfd} can be performed, i.e,
\begin{eqnarray}
\rho_{n,h}^{simp}=(\bar{U}_h^{n+1,simp}-P_h \bar{u}_{n+1}^{simp})+P_h \rho_n^{simp},
\label{rhonhs}
\end{eqnarray}
where $\bar{U}_h^{n+1,simp}$ is that corresponding to (\ref{Uhn}) with $U_h^n$ substituted by $P_h u(t_n)$, $p=1$ and the simplified value $\partial \bar{F}_{n,i}=\partial f(t_n,u(t_n))$. On the other hand,
$$\bar{u}_{n+1}^{simp}=\phi_{0,1,u(t_n),u(t_n)}(k)+k \sum_{i=1}^s \sum_{l=1}^s \mu_{i,l} \phi_{l,0,\bar{F}_{n,i},\bar{F}_{n,i}^{u,1,simp}}(k),$$
where $\bar{F}_{n,i}^{u,1,simp}=f(t_n,u(t_n))$, and then $\rho_n^{simp}=\bar{u}_{n+1}^{simp}-u(t_{n+1})$. (We notice that, for $p=1$, there is no difference in the required boundaries for the stages and therefore $\bar{F}_{n,i}$ is the same which would correspond to $\bar{u}_{n+1}$.)

Then, in (\ref{rhonhs}), the first parenthesis can be bounded in the same way than in the proof of Theorem \ref{thlocfd} and, in order to bound the second term, it suffices to consider that
$$\rho_n^{simp}=(\bar{u}_{n+1}^{simp}-\bar{u}_{n+1})+\rho_n,$$
that $\rho_n$ can be bounded as in Theorem \ref{thloc} with $p=1$ and that, using (\ref{phijr}),
\begin{eqnarray}
\bar{u}_{n+1}-\bar{u}_{n+1}^{simp} &=& k \sum_{i=1}^s \sum_{l=1}^s \mu_{i,l} \big[ \frac{1}{l!}[f(t_n+c_i k, \bar{K}_{n,i})-f(t_n,u(t_n))] \nonumber \\
&&\hspace{2cm}-\varphi_l(k A_0)[f(t_n+c_i k,\bar{K}_{n,i})-f(t_n,u(t_n))] \nonumber \\
&&\hspace{2cm}+k \varphi_{l+1}(k A_0) A[ f(t_n+c_i k, \bar{K}_{n,i})-f(t_n,u(t_n))] \big]. \label{difubu}
\end{eqnarray}
Then, by (\ref{bkni}), $\|\bar{K}_{n,i}-u(t_n)\|=O(k)$ and therefore, because of the assumed regularity of $\phi$, $h$ and (A7), the term in small brackets in (\ref{difubu}) is $O(k)$ and $A$ times that bracket is at least bounded. From this, it is clear that $\|\bar{u}_{n+1}-\bar{u}_{n+1}^{simp}\|=O(k^2)$, so that the first part of the theorem follows.

As for the second part of the theorem when $p=1$, it suffices to follow the same arguments as in the proof of Theorem \ref{thlocfd} again, taking into account that $A_0^{-1} \rho_n^{simp}=O(k^3)$ because, using (\ref{recurf}),
\begin{eqnarray}
A_0^{-1}(\bar{u}_{n+1}-\bar{u}_{n+1}^{simp}) &=& k \sum_{i=1}^s \sum_{l=1}^s \mu_{i,l} \bigg[ -k \varphi_{l+1}(k A_0)[f(t_n+c_i k, \bar{K}_{n,i})-f(t_n,u(t_n))] \nonumber \\
&&\hspace{1cm}+k^2 \varphi_{l+2}(k A_0) A[ f(t_n+c_i k, \bar{K}_{n,i})-f(t_n,u(t_n))] \nonumber \\
&&\hspace{1cm}+\frac{k}{(l+1)!} A_0^{-1} A[ f(t_n+c_i k, \bar{K}_{n,i})-f(t_n,u(t_n))] \bigg]. \nonumber
\end{eqnarray}
It suffices to notice that $\bar{K}_{n,i}-u(t_n)\in D(A_0)$ and because of (A8) and (A7),
\begin{eqnarray}
\lefteqn{\|A_0^{-1}A[ f(t_n+c_i k, \bar{K}_{n,i})-f(t_n,u(t_n))]\|} \nonumber \\
&&=\| \int_0^1 A_0^{-1}A[\phi'(\lambda \bar{K}_{n,i}+(1-\lambda)u(t_n))(\bar{K}_{n,i}-u(t_n))] d\lambda+c_i k A_0^{-1} A \dot{h}(t^*) \|=O(k),
\nonumber
\end{eqnarray}
from what the result follows.

As for $p=2$, by using again (\ref{rhonhs}), the key point is to bound appropriately $\rho_n^{simp}$ under the first new set of hypotheses and $A_0^{-1}\rho_n^{simp}$ under the second set. However, now the stages are calculated differently than in (\ref{Kbord}) since the approximation $\partial \bar{F}_{n,i}\approx \partial \bar{F}_{n,i}^{K,2,simp}$ at each stage with $\bar{F}_{n,i}^{K,2,simp}=f(t_n,u(t_n))$ must be used. More precisely,
\begin{eqnarray}
\bar{K}_{n,i}^{simp}&=&\phi_{0,1,u(t_n),u(t_n)}(c_i k)\nonumber \\
&&+k \sum_{j=1}^{i-1} \sum_{l,r=1}^s \lambda_{i,j,l,r} \phi_{l,0,\bar{F}_{n,j}^{simp},\bar{F}_{n,j}^{K,2,simp}}(c_r k), \quad i=1,\dots,s,
\label{bknis}
\end{eqnarray}
with $\bar{F}_{n,j}^{simp}=f(t_n+c_j k,\bar{K}_{n,j}^{simp})$.  On the other hand, as for the definition of $\bar{u}_{n+1}^{simp}$, the approximation for the boundaries in Table \ref{t2} must be considered, denoting
\begin{eqnarray}
\bar{F}_{n,i}^{u,2,simp}&=&f \big(t_n+c_i k,u(t_n) +k[c_i A u(t_n)+\sum_{j,l,r} \frac{\lambda_{i,j,l,r}}{l!} f(t_n,u(t_n))]\big), \nonumber \\
 \bar{F}_{n,i}^{u,2,A,simp}&=&f(t_n,u(t_n)),\nonumber
 \end{eqnarray}
it follows that
\begin{eqnarray}
\bar{u}_{n+1}-\bar{u}_{n+1}^{simp}=k \sum_{i,l} \mu_{i,l} \Delta_{i,l}(k),
\label{bus}
\end{eqnarray}
where
\begin{eqnarray}
\Delta_{i,l}'(\tau)&=&(A-\frac{l}{\tau}I) \Delta_{i,l}(\tau)+\frac{1}{(l-1)!\tau}[\bar{F}_{n,i}-\bar{F}_{n,i}^{simp}], \nonumber \\
\Delta_{i,l}(0)&=&\frac{1}{l!}[\bar{F}_{n,i}-\bar{F}_{n,i}^{simp}],  \nonumber \\
\partial \Delta_{i,l}(\tau)&=&\partial \bigg[\frac{1}{l!}[\bar{F}_{n,i}-\bar{F}_{n,i}^{u,2,simp}]+\frac{\tau}{(l+1)!}A[ \bar{F}_{n,i}- \bar{F}_{n,i}^{u,2,A,simp}] \bigg]. \label{bordD}
\end{eqnarray}
From here, denoting by $W_{i,l}(\tau)$ to the term in big brackets in (\ref{bordD}),
\begin{eqnarray}
(\Delta_{i,l}-W_{i,l})'(\tau)&=&(A-\frac{l}{\tau}I) (\Delta_{i,l}-W_{i,l}) (\tau)+\frac{1}{(l-1)!\tau}[\bar{F}_{n,i}^{u,2,simp}-\bar{F}_{n,i}^{simp}] \nonumber \\
&&+\frac{1}{l!}A[\bar{F}_{n,i}^{u,A,2,simp}-\bar{F}_{n,i}^{u,2,simp}]+\frac{\tau}{(l+1)!}A^2[\bar{F}_{n,i}-\bar{F}_{n,i}^{u,2,simp}], \nonumber \\
(\Delta_{i,l}-W_{i,l})(0)&=&\frac{1}{l!}[\bar{F}_{n,i}^{u,2,simp}-\bar{F}_{n,i}^{simp}],  \nonumber \\
\partial (\Delta_{i,l}-W_{i,l})(\tau)&=&0, \nonumber
\end{eqnarray}
so that
\begin{eqnarray}
(\Delta_{i,l}-W_{i,l})(k)&=&\varphi_l(k A_0)[\bar{F}_{n,i}^{u,2,simp}-\bar{F}_{n,i}^{simp}]+k\varphi_{l+1}(k A_0) A[\bar{F}_{n,i}^{u,2,A,simp}-\bar{F}_{n,i}^{u,2,simp}] \nonumber \\
&&+k^2 \varphi_{l+2}(k A_0) A^2[\bar{F}_{n,i}-\bar{F}_{n,i}^{u,2,A,simp}]. \label{DW}
\end{eqnarray}
Now, from (\ref{bknis}), (\ref{phi0r}), (\ref{phijr}), the assumed regularity and (A7), $\|\bar{F}_{n,i}^{u,2,simp}-\bar{F}_{n,i}^{simp}\|=O(k^2)$. On the other hand, it is also clear that $\|A[\bar{F}_{n,i}^{u,A,2,simp}-\bar{F}_{n,i}^{u,2,simp}]\|=O(k)$. Therefore, from (\ref{DW}), $\|(\Delta_{i,l}-W_{i,l})(k)\|=O(k^2)$. Using also that, from (\ref{bkni}),  $\bar{F}_{n,i}-\bar{F}_{n,i}^{u,2,s}=O(k^2)$ and through (A7), $\|A[\bar{F}_{n,i}-\bar{F}_{n,i}^{u,2,A,simp}]\|=O(k)$, it follows that $W_{i,l}(k)$ is also $O(k^2)$.  Then,
$$\|\Delta_{i,l}(k)\|=\|(\Delta_{i,l}-W_{i,l})(k)+W_{i,l}(k)\|=O(k^2),$$
and from (\ref{bus}) and the first part of Theorem \ref{thloc} for $p=2$, $\rho_n^{simp}=O(k^3)$, which proves the first part of the present theorem.

In order to prove the second part of the theorem for $p=2$, let us first notice that, using (\ref{phijr}),
\begin{eqnarray}
\bar{K}_{n,i}-\bar{K}_{n,i}^{simp} &=& k \sum_{j=1}^{i-1} \sum_{l,r=1}^s \lambda_{i,j,l,r} \big[ \frac{1}{l!}[\bar{F}_{n,j}-\bar{F}_{n,j}^{K,2,simp}]-\varphi_l(c_r k A_0)[\bar{F}_{n,j}^{simp}-\bar{F}_{n,j}^{K,2,simp}] \nonumber \\
&&\hspace{3cm}+c_r k \varphi_{l+1}(c_r k A_0) A[ \bar{F}_{n,j}-\bar{F}_{n,j}^{K,2,simp}] \big]. \nonumber
\end{eqnarray}
Inductively on the stages, it is therefore proved that $\|\bar{K}_{n,i}-\bar{K}_{n,i}^{simp}\|=O(k^2)$, from what $\|\bar{F}_{n,i}-\bar{F}_{n,i}^{simp}\|=O(k^2)$. But, furthermore, using (\ref{recurf}),
\begin{eqnarray}
A_0^{-1}(\bar{K}_{n,i}-\bar{K}_{n,i}^{simp}) &=& k \sum_{j=1}^{i-1} \sum_{l,r=1}^s \lambda_{i,j,l,r} \bigg[ \frac{1}{l!}A_0^{-1}[\bar{F}_{n,j}-\bar{F}_{n,j}^{simp}]-k \varphi_{l+1}(c_r k A_0)
[\bar{F}_{n,j}^{simp}-\bar{F}_{n,j}^{K,2,simp}]
\nonumber \\
&&\hspace{0.5cm}+(c_r k)^2 \varphi_{l+2}(c_r k A_0) A[ \bar{F}_{n,j}-\bar{F}_{n,j}^{K,2,simp}]
+ \frac{c_r k}{(l+1)!}A_0^{-1} A [\bar{F}_{n,j}-\bar{F}_{n,j}^{K,2,simp}]\bigg], \nonumber
\end{eqnarray}
which implies, using (A7) and (A8), that $\|A_0^{-1}(\bar{K}_{n,i}-\bar{K}_{n,i}^{simp})\|=O(k^3)$.
Now, we notice that, from (\ref{DW}) and the definition of $W_{i,l}$,
\begin{eqnarray}
\Delta_{i,l}(k)&=&(\varphi_l(k A_0)-\frac{1}{l!}I)[\bar{F}_{n,i}^{u,2,simp}-\bar{F}_{n,i}^{simp}]+\frac{1}{l!}[\bar{F}_{n,i}-\bar{F}_{n,i}^{simp}] \nonumber \\
&&+k \varphi_{l+1}(k A_0) A [\bar{F}_{n,i}^{u,2,A,simp}-\bar{F}_{n,i}^{u,2,simp}] +\frac{k}{(l+1)!}A [\bar{F}_{n,i}-\bar{F}_{n,i}^{u,2,simp}] \nonumber \\
&&+k^2 \varphi_{l+2}(k A_0)A^2 [\bar{F}_{n,i}-\bar{F}_{n,i}^{u,2,A,simp}],
\nonumber
\end{eqnarray}
from what, applying $A_0^{-1}$ and then (\ref{recurf}),
\begin{eqnarray}
A_0^{-1} \Delta_{i,l}(k)&=& k \varphi_{l+1}(k A_0)[\bar{F}_{n,i}^{u,2,simp}-\bar{F}_{n,i}^{simp}]+\frac{1}{l!}A_0^{-1}[\bar{F}_{n,i}-\bar{F}_{n,i}^{simp}] \nonumber \\
&&+k^2 \varphi_{l+2}(k A_0) A [\bar{F}_{n,i}^{u,2,A,simp}-\bar{F}_{n,i}^{u,2,simp}] +\frac{k}{(l+1)!}A_0^{-1} A [\bar{F}_{n,i}-\bar{F}_{n,i}^{u,2,simp}] \nonumber \\
&&+k^2 \varphi_{l+2}(k A_0) A_0^{-1} A^2 [\bar{F}_{n,i}-\bar{F}_{n,i}^{u,2,A,simp}].
\nonumber
\end{eqnarray}
Now, the first and third terms are $O(k^3)$ because of previously justified bounds. The second term is $O(k^3)$ using (\ref{cspc}), (A7) and the fact that $\|A_0^{-1}(\bar{K}_{n,i}-\bar{K}_{n,i}^{simp})\|=O(k^3)$. The fourth term is $O(k^3)$ because, using (A8), $\|A_0^{-1} A [\bar{F}_{n,i}-\bar{F}_{n,i}^{u,2,simp}]\|=O(k^2)$. The fifth term is $O(k^3)$ because $\|A(\bar{K}_{n,i}-u(t_n))\|=O(k)$ as it is deduced from the proof of Theorem \ref{thloc} and then, applying (A8), $\|A_0^{-1}A^2[\bar{F}_{n,i}-\bar{F}_{n,i}^{u,2,A,simp}]\|=O(k)$. Summing up, $\|A_0^{-1}\Delta_{i,l}(k)\|=O(k^3)$, from what $A_0^{-1}(\bar{u}_{n+1}-\bar{u}_{n+1}^{simp})$ and $A_0^{-1}\rho_n^{simp}$ are $O(k^4)$.

As for $p=3$, the proof follows the same lines of the case $p=2$ with the difference that now,
\begin{eqnarray}
\bar{K}_{n,i}^{simp}=\phi_{0,2,u(t_n),u(t_n)}(c_i k)+k\sum_{j=1}^{i-1} \sum_{l,r=1}^s \lambda_{i,j,l,r} \tilde{\phi}_{l,j}(c_r k),
\label{Knis}
\end{eqnarray}
where
\begin{eqnarray}
\tilde{\phi}_{l,j}'(\tau)&=&(A-\frac{l}{\tau})\tilde{\phi}_{l,j}(\tau)+\frac{1}{(l-1)!\tau}\bar{F}_{n,j}^{simp}, \nonumber \\
\tilde{\phi}_{l,j}(0)&=&\frac{1}{l!}\bar{F}_{n,j}^{simp}, \nonumber \\
\partial \tilde{\phi}_{l,j}(\tau)&=& \partial [ \frac{1}{l!}\bar{F}_{n,j}^{K,3,simp}+\frac{\tau}{(l+1)!}A \bar{F}_{n,j}^{K,3,A,simp}], \nonumber
\end{eqnarray}
with
$\bar{F}_{n,j}^{K,3,simp}$ and $\bar{F}_{n,j}^{K,3,A,simp}$ the corresponding functions taken from Table \ref{t1}:
\begin{eqnarray}
\bar{F}_{n,j}^{K,3,simp}&=&f \big(t_n+c_j k,u(t_n) +k[c_j A u(t_n)+\sum_{jj=1}^{j-1}\sum_{l,r} \frac{\lambda_{j,jj,l,r}}{l!} f(t_n,u(t_n))]\big), \label{fnjk3} \\
\bar{F}_{n,j}^{K,3,A,simp}&=&f(t_n,u(t_n)). \nonumber
\end{eqnarray}
From this, considering the notation
\begin{eqnarray}
T_{l,j}(\tau)=\frac{1}{l!}\bar{F}_{n,j}^{K,3,simp}+\frac{\tau}{(l+1)!}A \bar{F}_{n,j}^{K,3,A,simp},
\nonumber
\end{eqnarray}
the following is satisfied
\begin{eqnarray}
(\tilde{\phi}_{l,j}-T_{l,j})'(\tau)&=&(A-\frac{l}{\tau})(\tilde{\phi}_{l,j}-T_{l,j})(\tau)+\frac{1}{(l-1)!\tau}[\tilde{F}_{n,j}^{simp}-\bar{F}_{n,j}^{K,3,simp}] \nonumber \\
&&+\frac{1}{l!}A[\bar{F}_{n,j}^{K,3,simp}-\bar{F}_{n,j}^{K,3,A,simp}]+\frac{\tau}{(l+1)!}A^2 \bar{F}_{n,j}^{K,3,A,simp}, \nonumber \\
(\tilde{\phi}_{l,j}-T_{l,j})(0)&=&\frac{1}{l!}[\tilde{F}_{n,j}^{simp}-\tilde{F}_{n,j}^{K,3,simp}], \nonumber \\
\partial (\tilde{\phi}_{l,j}-T_{l,j})(\tau)&=& 0. \nonumber
\end{eqnarray}
Solving this,
\begin{eqnarray}
(\tilde{\phi}_{l,j}-T_{l,j})(c_r k)&=& \varphi_l(c_r k A_0)[\bar{F}_{n,j}^{simp}-\bar{F}_{n,j}^{K,3,simp}]+c_r k \varphi_{l+1}(c_r k A_0) A [\bar{F}_{n,j}^{K,3,simp}-\bar{F}_{n,j}^{K,3,A,simp}] \nonumber \\
&&+(c_r k)^2 \varphi_{l+2}(\tau A_0) A^2 \tilde{F}_{n,j}^{K,3,A,simp} \nonumber \\
&=&O(k^2), \nonumber
\end{eqnarray}
where, for the last equality, an induction argument over the stages has been used, (\ref{Knis}), (\ref{fnjk3}), as well as the assumed regularity and (A7). From this,
\begin{eqnarray}
\bar{K}_{n,i}^{simp}&=&u(t_n)+c_i k A u(t_n)+\frac{(c_i k)^2}{2}A^2 u(t_n)\nonumber \\
&&+k \sum_{j=1}^{i-1} \sum_{l,r=1}^s \big[ \frac{1}{l!}\bar{F}_{n,j}^{K,3,simp}+\frac{c_r k}{(l+1)!}A \bar{F}_{n,j}^{K,3,A,simp} \big]+O(k^3). \label{kbsimp4}
\end{eqnarray}
On the other hand,
\begin{eqnarray}
\bar{u}_{n+1}-\bar{u}_{n+1}^{simp}=k \sum_{i,l=1}^s \mu_{i,l} \tilde{\Delta}_{i,l}(k), \label{ubn3}
\end{eqnarray}
where
\begin{eqnarray}
\tilde{\Delta}_{i,l}'(\tau)&=&(A-\frac{l}{\tau})\Delta_{i,l}(\tau)+\frac{1}{(l-1)! \tau}[\bar{F}_{n,i}-\bar{F}_{n,i}^{simp}], \nonumber \\
\tilde{\Delta}_{i,l}(0)&=&\frac{1}{l!}[\bar{F}_{n,i}-\bar{F}_{n,i}^{simp}], \nonumber \\
\partial \tilde{\Delta}_{i,l}(\tau)&=&\partial \bigg[\frac{1}{l!}[\bar{F}_{n,i}-\bar{F}_{n,i}^{u,3,simp}]+\frac{\tau}{(l+1)!}A [\bar{F}_{n,i}-\bar{F}_{n,i}^{u,3,A,simp}] 
+\frac{\tau^2}{(l+2)!}A^2 [\bar{F}_{n,i}-\bar{F}_{n,i}^{u,3,A^2,simp}]\bigg], \nonumber
\end{eqnarray}
where $\bar{F}_{n,i}^{u,3,simp}$, $\bar{F}_{n,i}^{u,3,A,simp}$ and $\bar{F}_{n,i}^{u,3,A^2,simp}$ are again the corresponding expressions in Table \ref{t2}. Considering then the notation
\begin{eqnarray}
\tilde{W}_{i,l}(\tau)&=& \frac{1}{l!}[\bar{F}_{n,i}-\bar{F}_{n,i}^{u,3,simp}]+\frac{\tau}{(l+1)!}A [\bar{F}_{n,i}-\bar{F}_{n,i}^{u,3,A,simp}] \nonumber \\
 &&+\frac{\tau^2}{(l+2)!}A^2 [\bar{F}_{n,i}-\bar{F}_{n,i}^{u,3,A^2,simp}], \label{tWil}
\end{eqnarray}
it follows that
\begin{eqnarray}
(\tilde{\Delta}_{i,l}-\tilde{W}_{i,l})'(\tau)&=&(A-\frac{l}{\tau})(\tilde{\Delta}_{i,l}-\tilde{W}_{i,l})(\tau)+\frac{1}{l!}A[\bar{F}_{n,i}^{u,3,A,simp}-\bar{F}_{n,i}^{u,3,simp}] \nonumber \\
&&+\frac{\tau}{(l+1)!}A^2[\bar{F}_{n,i}^{u,3,A^2,simp}-\bar{F}_{n,i}^{u,3,A,simp}]+\frac{\tau^2}{(l+2)!}A^3[\bar{F}_{n,i}-\bar{F}_{n,i}^{u,3,A^2,simp}] \nonumber \\
&&+\frac{1}{(l-1)! \tau}[\bar{F}_{n,i}^{u,3,simp}-\bar{F}_{n,i}^{simp}], \nonumber \\
(\tilde{\Delta}_{i,l}-\tilde{W}_{i,l})(0)&=&\frac{1}{l!}[\bar{F}_{n,i}^{u,3,simp}-\bar{F}_{n,i}^{simp}], \nonumber \\
\partial (\tilde{\Delta}_{i,l}-\tilde{W}_{i,l})(\tau)&=&0. \nonumber
\end{eqnarray}
From this,
\begin{eqnarray}
(\tilde{\Delta}_{i,l}-\tilde{W}_{i,l})(k)&=&\varphi_l(k A_0)[\bar{F}_{n,i}^{u,3,simp}-\bar{F}_{n,i}^{simp}]+k \varphi_{l+1}(k A_0)A [\bar{F}_{n,i}^{u,3,A,simp}-\bar{F}_{n,i}^{u,3,simp}] \nonumber \\
&&+k^2 \varphi_{l+2}(k A_0)A^2 [\bar{F}_{n,i}^{u,3,A^2,simp}-\bar{F}_{n,i}^{u,3,A,simp}] \nonumber \\
&&+k^3 \varphi_{l+3}(k A_0) A^3 [\bar{F}_{n,i}-\bar{F}_{n,i}^{u,3,A^2,simp}].
\label{DWtilde}
\end{eqnarray}
Considering (\ref{kbsimp4}), (A7) and the assumed regularity, $\|\bar{F}_{n,i}^{u,3,simp}-\bar{F}_{n,i}^{simp}\|=O(k^3)$. In a similar way,  the rest of the terms which turn up here, as well as those in (\ref{tWil}) for $\tau=k$ are $O(k^3)$, so that $\tilde{\Delta}_{i,l}(k)=O(k^3)$. Therefore, using (\ref{ubn3}), $\bar{u}_{n+1}-\bar{u}_{n+1}^{simp}=O(k^4)$, from what, using Theorem \ref{thloc} with $p=3$, $\rho_n^{simp}=O(k^4)$ and, finally, through the same proof of Theorem \ref{thlocfd} with $p=3$, the first result of the present theorem follows.

As for the second part of the theorem, we firstly notice that
\begin{eqnarray}
A_0^{-1}(\bar{K}_{n,i}-\bar{K}_{n,i}^{simp})=k \sum_{j=1}^{i-1} \sum_{l,r=1}^s \lambda_{i,j,l,r} A_0^{-1} \Delta_{j,l}(c_r k), \label{invk3}
\end{eqnarray}
where
$\Delta_{j,l}$ is the function in the proof of Theorem \ref{thlocs}. From there, $A_0^{-1} \Delta_{j,l}(c_r k)=O(k^3)$, which implies that $A_0^{-1}(\bar{K}_{n,i}-\bar{K}_{n,i}^{simp})=O(k^4)$. Then,
\begin{eqnarray}
A_0^{-1}(\bar{u}_{n+1}-\bar{u}_{n+1}^{simp})=k \sum_{i,l} \mu_{i,l} A_0^{-1} \tilde{\Delta}_{i,l}(k),
\label{invu3}
\end{eqnarray}
and from (\ref{DWtilde}) and (\ref{tWil}),
\begin{eqnarray}
A_0^{-1} \tilde{\Delta}_{i,l}(k)&=&\frac{1}{l!}A_0^{-1}(\bar{F}_{n,i}-\bar{F}_{n,i}^{simp})+k \varphi_{l+1}(k A_0)[\bar{F}_{n,i}^{u,3,simp}-\bar{F}_{n,i}^{simp}]
\nonumber \\
&&+k^2 \varphi_{l+2}(k A_0) A [\bar{F}_{n,i}^{u,3,A,simp}-\bar{F}_{n,i}^{u,3,simp}]+\frac{k}{(l+1)!}A_0^{-1} A[\bar{F}_{n,i}-\bar{F}_{n,i}^{u,3,simp}] \nonumber \\
&&+k^3 \varphi_{l+3}(k A_0) A^2 [\bar{F}_{n,i}^{u,3,A^2,simp}-\bar{F}_{n,i}^{u,3,A,simp}]+\frac{k^2}{(l+2)!}A_0^{-1} A^2[\bar{F}_{n,i}-\bar{F}_{n,i}^{u,3,A,simp}] \nonumber \\
&&+k^3 \varphi_{l+3}(k A_0) A_0^{-1} A^3 [\bar{F}_{n,i}-\bar{F}_{n,i}^{u,3,A^2,simp}]. \nonumber
\end{eqnarray}
Now, the first term in the right hand side is $O(k^4)$ using (\ref{cspc}) and (\ref{invk3}). The second, third and fifth terms are also $O(k^4)$ using the assumed regularity and (A7). Finally, the fourth, sixth and seventh terms are $O(k^4)$ by using (A8) and the fact that $\| \bar{K}_{n,i}-\bar{K}_{n,i}^{u,3,simp}\|_h=O(k^3)$, $\| A(\bar{K}_{n,i}-\bar{K}_{n,i}^{u,3,A,simp})\|_h=O(k^2)$ and $\| A^2(\bar{K}_{n,i}-\bar{K}_{n,i}^{u,3,A^2,simp})\|_h=O(k)$, considering the proof of Theorem \ref{thloc}.
Therefore, from (\ref{invu3}), $A_0^{-1}(\bar{u}_{n+1}-\bar{u}_{n+1}^{simp})=O(k^5)$, from what $A_0^{-1}\rho_n^{simp}=O(k^5)$, and the result follows considering the same proof as that of the second part of Theorem \ref{thlocfd} with $p=3$.

\subsection{Proof of Theorem \ref{thglob}}

In the following $\rho_n$ and $\rho_{n,h}$ in fact correspond to $\rho_n^{sim}$ and $\rho_{n,h}^{sim}$.
By definition,
\begin{eqnarray}
e_{n+1,h}=(U_h^{n+1}-\bar{U}_h^{n+1})+\rho_{n,h}.
\label{decompe}
\end{eqnarray}
Then, for $p=1$, using Table \ref{t3},
\begin{eqnarray}
U_h^{n+1}-\bar{U}_h^{n+1}&=&e^{ k A_{h,0}}(U_h^n-P_h u(t_n))-k \varphi_1(k A_{h,0}) D_h O(e_{n,h})
+k^2 \varphi_2(k A_{h,0}) C_h O(e_{n,h}) \nonumber \\
&&+k \sum_{i,l=1}^s \mu_{i,l}\big[\varphi_l(k A_{h,0})[F_{n,i,h}-\bar{F}_{n,i,h}]+k \varphi_{l+1}(k A_{h,0})C_h  O(e_{n,h})\big], \nonumber
\end{eqnarray}
where the term $C_h O(e_{n,h})$ is just non-zero for Robin/Neumann boundary conditions. On the other hand, using (\ref{Knih1}),
\begin{eqnarray}
K_{n,i,h}-\bar{K}_{n,i,h}&=&e^{c_i k A_{h,0}}(U_h^n-P_h u(t_n))+k \sum_{j=1}^{i-1} \sum_{l,r=1}^s \lambda_{i,j,l,r} \varphi_l(c_r k A_{h,0})[F_{n,j,h}-\bar{F}_{n,j,h}],
\nonumber
\end{eqnarray}
and so, inductively, because of (H3), $F_{n,i,h}-\bar{F}_{n,i,h}=O(e_{n,h})$.
Inserting this in (\ref{decompe}), using also (\ref{recurf}), (H2b) and (H2c),
$$
e_{n+1,h}=e^{k A_{h,0}}e_{n,h}+O(k e_{n,h})+\rho_{n+1,h}.
$$
Then, applying Gronwall lemma, under the first set of hypotheses of Theorem \ref{thlocfd}, using the bound for $\rho_{n,h}$ and also Theorem \ref{thlocs},
the result for $p=1$ follows. Under the second set of hypotheses, the result for $p=1$ can be obtained using a summation-by-parts argument, which can be applied in a similar way to what was done for exponential splitting methods in \cite{acrnl} taking (\ref{parabol}), the regularity (\ref{regflesp1}) and the bound for $A_{h,0}^{-1}\rho_{n,h}$ into account.

For $p=2$, we will separate the cases corresponding to Dirichlet and to R/N boundary conditions.

With Dirichlet boundary conditions, considering Table \ref{t4},
\begin{eqnarray}
\lefteqn{\hspace{-1cm}U_h^{n+1}-\bar{U}_h^{n+1}=e^{ k A_{h,0}}(U_h^n-P_h u(t_n))
-k^2 \varphi_2(k A_{h,0}) D_h O(\nu_h+\frac{e_{n,h}}{h^\gamma})}\nonumber \\
&&+k^3 \varphi_3(k A_{h,0}) C_h O(\nu_h+\frac{e_{n,h}}{h^\gamma}) \nonumber \\
&&+k \sum_{i,l=1}^s \mu_{i,l}\bigg[\varphi_l(k A_{h,0})[F_{n,i,h}-\bar{F}_{n,i,h}]-k \varphi_{l+1}(k A_{h,0}) D_h O(\nu_h+\frac{e_{n,h}}{h^\gamma})
\nonumber \\
&& \hspace{2cm}+k^2 \varphi_{l+2}(k A_{h,0})C_h  O(\nu_h+\frac{e_{n,h}}{h^\gamma})\bigg], \nonumber
\end{eqnarray}
where $F_{n,i,h}-\bar{F}_{n,i,h}=O(e_{n,h})$ with the same argument of the previous theorem. Then, using (H2b) and (H2c),
$$
e_{n+1,h}=e^{k A_{h,0}}e_{n,h}+O(k^2 \nu_h+\frac{k^2}{h^\gamma}e_{n,h}+k e_{n,h})+\rho_{n+1,h},
$$
from what, using also (\ref{cfl}) and the first set of hypotheses of Theorem \ref{thlocs},
$$
e_{n+1,h}=e^{k A_{h,0}}e_{n,h}+O(k^2 \nu_h+k e_{n,h}+k^3+k\varepsilon_h),
$$
and the result follows by a discrete Gronwall lemma. Under the second set of hypotheses, the summation-by-parts argument can be applied as in the previous theorem using the respective results of Theorems \ref{thlocfd} and \ref{thlocs}.

On the other hand, with R/N boundary conditions and $p=2$, using again Table \ref{t4},
\begin{eqnarray}
U_h^{n+1}-\bar{U}_h^{n+1}
&=&e^{k A_{h,0}}(U_h^n-P_h u(t_n))-k \varphi_1(k A_{h,0}) D_h O(e_{n,h}) \nonumber \\
&&+k^2 \varphi_2(k A_{h,0})[C_h O(e_{n,h})-D_h O(\mu_{k,1}+\frac{e_{n,h}}{k}+\nu_h+\frac{e_{n,h}}{h^\gamma})]\nonumber \\
&&+k^3 \varphi_3(k A_{h,0})C_h O(\mu_{k,1}+\frac{e_{n,h}}{h^\gamma}+\nu_h+\frac{e_{n,h}}{k}) \nonumber \\
&&+k \sum_{i,l=1}^s \mu_{i,l} \bigg[\varphi_l(k A_{h,0})[F_{n,i,h}-\bar{F}_{n,i,h}] \nonumber \\
&&\hspace{2cm}+k \varphi_{l+1}(k A_{h,0})[C_h O(e_{n,h}+k \mu_{k,1})-D_h O(\mu_{k,1}+\frac{e_{n,h}}{k}+\nu_h+\frac{e_{n,h}}{h^\gamma})] \nonumber \\
&&\hspace{2cm}+k^2 \varphi_{l+2}(k A_{h,0})C_h  O(\mu_{k,1}+\frac{e_{n,h}}{k}+\nu_h+\frac{e_{n,h}}{h^\gamma})\bigg], \nonumber
\end{eqnarray}
where
\begin{eqnarray}
K_{n,i,h}-\bar{K}_{n,i,h}&=&e^{c_i k A_{h,0}}(U_h^n-P_h u(t_n))+(c_i k)^2 \varphi_2(c_i k A_{h,0})C_h O(e_{n,h}) \nonumber \\
&&+k \sum_{j=1}^{i-1} \sum_{l,r=1}^s \lambda_{i,j,l,r} \big[ \varphi_l(c_r k A_{h,0})[F_{n,j,h}-\bar{F}_{n,j,h}]+ c_r k \varphi_{l+1}(c_r k A_{h,0})C_h O(e_{n,h})\big],
\nonumber
\end{eqnarray}
and so again, inductively, because of (H2b) and (H3), $F_{n,i,h}-\bar{F}_{n,i,h}=O(e_{n,h})$. Then,
$$
e_{n+1,h}=e^{k A_{h,0}}e_{n,h}+O(k e_{n,h}+k^2 \mu_{k,1}+k^2 \nu_h+\frac{k^2}{h^\gamma}e_{n,h})+\rho_{n+1,h},
$$
and arguing as before, the result follows.

For $p=3$, the theorem follows as before but taking now into account that, with Dirichlet boundary conditions,
using  Table \ref{t5},
\begin{eqnarray}
U_h^{n+1}-\bar{U}_h^{n+1}&=&e^{k A_{h,0}}(U_h^n-P_h u(t_n))-k^2 \varphi_2(k A_{h,0})D_h O(\nu_h+\frac{e_{n,h}}{h^\gamma}) \nonumber \\
&&+k^3 \varphi_3(k A_{h,0})[C_h O(\nu_h+\frac{e_{n,h}}{h^\gamma})-D_h O(\nu_h+\frac{e_{n,h}}{k h^\gamma}+\frac{\mu_{k,1}}{h^\gamma})]\nonumber \\
&&+k^4 \varphi_4(k A_{h,0})C_h O(\nu_h+\frac{e_{n,h}}{k h^\gamma}+\frac{\mu_{k,1}}{h^\gamma}) \nonumber \\
&&+k \sum_{i,l=1}^s \mu_{i,l} \bigg[\varphi_l(k A_{h,0})[F_{n,i,h}-\bar{F}_{n,i,h}] \nonumber \\
&&\hspace{2cm}+k \varphi_{l+1}(k A_{h,0})[C_h O(k^2 \nu_h+\frac{k^2}{h^\gamma}e_{n,h})-D_h O(\nu_h+\frac{e_{n,h}}{h^\gamma}+\frac{k \mu_{k,1}}{h^\gamma})] \nonumber \\
&&\hspace{2cm}+k^2 \varphi_{l+2}(k A_{h,0})[C_h  O(\nu_h+\frac{e_{n,h}}{h^\gamma}+\frac{k \mu_{k,1}}{h^\gamma})-D_h  O(\nu_h+\frac{e_{n,h}}{k h^\gamma}+\frac{\mu_{k,1}}{h^\gamma})] \nonumber \\
&&\hspace{2cm}+k^3 \varphi_{l+3}(k A_{h,0})C_h  O(\nu_h+\frac{e_{n,h}}{k h^\gamma}+\frac{\mu_{k,1}}{h^\gamma})\bigg], \nonumber
\end{eqnarray}
where
\begin{eqnarray}
K_{n,i,h}-\bar{K}_{n,i,h}&=&e^{c_i k A_{h,0}}(U_h^n-P_h u(t_n))-(c_i k)^2 \varphi_2(k A_{h,0})D_h O(\nu_h+\frac{e_{n,h}}{h^\gamma}) \nonumber \\
&&+(c_i k)^3 \varphi_3(c_i k A_{h,0})C_h O(\nu_h+\frac{e_{n,h}}{h^\gamma}) \nonumber \\
&&+k \sum_{j=1}^{i-1} \sum_{l,r=1}^s \lambda_{i,j,l,r} \big[ \varphi_l(c_r k A_{h,0})[F_{n,j,h}-\bar{F}_{n,j,h}]-c_r k \varphi_{l+1}(c_r k A_{h,0}) D_h O(\nu_h+\frac{e_{n,h}}{h^\gamma})
\nonumber \\
&& \hspace{3cm}
+ (c_r k)^2 \varphi_{l+2}(c_r k A_{h,0})C_h O(\nu_h+\frac{e_{n,h}}{h^\gamma})\big],
\nonumber
\end{eqnarray}
and therefore $F_{n,i,h}-\bar{F}_{n,i,h}=O(e_{n,h}+k^2 \nu_h+\frac{k^2}{h^\gamma} e_{n,h})$. Using this,
$$
e_{n+1,h}=e^{k A_{h,0}}e_{n,h}+O(k e_{n,h}+k^2 \nu_h+\frac{k^2}{h^\gamma} e_{n,h}+k^3 \frac{\mu_{k,1}}{h^\gamma})+\rho_{n+1,h},
$$
from what arguing as in previous theorems, but using now the results in Theorem \ref{thlocs}, the result follows under the first and second set of hypotheses.

With R/N boundary conditions, using again Table \ref{t5},
{\footnotesize \begin{eqnarray}
\lefteqn{U_h^{n+1}-\bar{U}_h^{n+1}} \nonumber \\
&=&e^{k A_{h,0}}(U_h^n-P_h u(t_n))-k \varphi_1(k A_{h,0})D_h O(e_{n,h}) \nonumber \\
&&+k^2 \varphi_2(k A_{h,0})[C_h O(e_{n,h})-D_h O(\mu_{k,1}+\frac{e_{n,h}}{k}+\nu_h+\frac{e_{n,h}}{h^\gamma})]\nonumber \\
&&+k^3 \varphi_3(k A_{h,0})[C_h O(\mu_{k,1}+\frac{e_{n,h}}{k}+\nu_h+\frac{e_{n,h}}{h^\gamma})-D_h  O(\mu_{k,1}+\mu_{k,2}+\frac{e_{n,h}}{k^2}+\nu_h+\frac{e_{n,h}}{k h^\gamma})] \nonumber \\
&&+k^4 \varphi_4(k A_{h,0})C_h O(\mu_{k,1}+\mu_{k,2}+\frac{e_{n,h}}{k^2}+\nu_h+\frac{e_{n,h}}{k h^\gamma}) \nonumber \\
&&+k \sum_{i,l=1}^s \mu_{i,l} \big[\varphi_l(k A_{h,0})[F_{n,i,h}-\bar{F}_{n,i,h}] \nonumber \\
&&\hspace{1cm}+k \varphi_{l+1}(k A_{h,0})[C_h O(e_{n,h}+k \mu_{k,1}+k^2 \mu_{k,2}+k^2\nu_h+\frac{k^2}{h^\gamma}e_{n,h})-D_h O(\mu_{k,1}+k \mu_{k,2}+\frac{e_{n,h}}{k}+\nu_h+\frac{e_{n,h}}{h^\gamma})] \nonumber \\
&&\hspace{1cm}+k^2 \varphi_{l+2}(k A_{h,0})[C_h  O(\mu_{k,1}+k \mu_{k,2}+\frac{e_{n,h}}{k}+\nu_h+\frac{e_{n,h}}{h^\gamma})-D_h O(\mu_{k,1}+\mu_{k,2}+\frac{e_{n,h}}
{k^2}+\nu_h+\frac{e_{n,h}}{k h^\gamma})] \nonumber \\
&&\hspace{1cm}+k^3 \varphi_{l+3}(k A_{h,0})C_h  O(\mu_{k,1}+\mu_{k,2}+\frac{e_{n,h}}
{k^2}+\nu_h+\frac{e_{n,h}}{k h^\gamma})\big], \nonumber
\end{eqnarray}}
where
{\footnotesize \begin{eqnarray}
K_{n,i,h}-\bar{K}_{n,i,h}&=&e^{c_i k A_{h,0}}(U_h^n-P_h u(t_n))-c_i k \varphi_1(c_i k A_{h,0})D_h O(e_{n,h}) \nonumber \\
&&+(c_i k)^2 \varphi_2(c_i k A_{h,0})[C_h O(e_{n,h})-D_h O(\mu_{k,1}+\frac{e_{n,h}}{k}+\nu_h+\frac{e_{n,h}}{h^\gamma})] \nonumber \\
&&+(c_i k)^3 \varphi_3(c_i k A_{h,0})[C_h O(\mu_{k,1}+\frac{e_{n,h}}{k}+\nu_h+\frac{e_{n,h}}{h^\gamma})  \nonumber \\
&&+k \sum_{j=1}^{i-1} \sum_{l,r=1}^s \lambda_{i,j,l,r} \bigg[ \varphi_l(c_r k A_{h,0})[F_{n,j,h}-\bar{F}_{n,j,h}]
\nonumber \\
&& \hspace{1cm}+c_r k \varphi_{l+1}(c_r k A_{h,0})[C_h O(e_{n,h}+k \mu_{k,1})-D_h O(\mu_{k,1}+\frac{e_{n,h}}{k}+\nu_h+\frac{e_{n,h}}{h^\gamma})] \nonumber \\
&& \hspace{1cm}+ (c_r k)^2 \varphi_{l+2}(c_r k A_{h,0})C_h O(\mu_{k,1}+\frac{e_{n,h}}{k}+\nu_h+\frac{e_{n,h}}{h^\gamma})\bigg],
\nonumber
\end{eqnarray}}
and therefore $F_{n,i,h}-\bar{F}_{n,i,h}=O(e_{n,h}+k^2 \mu_{k,1}+k^2 \nu_h+\frac{k^2}{h^\gamma} e_{n,h})$. Using this,
$$
e_{n+1,h}=e^{k A_{h,0}}e_{n,h}+O(k e_{n,h}+k^2 \mu_{k,1}+k^2 \nu_h+\frac{k^2}{h^\gamma} e_{n,h}+k^3 \mu_{k,2})+\rho_{n+1,h},
$$
from what the result follows as before.

\end{document}